\documentclass{SLC}
\articlenumber{10}

\usepackage{enumitem}
\usepackage{tikz}
\usepackage[all]{xy}

\def\Z{\mathbb{Z}}
\def\N{\mathbb{N}}
\def\C{\mathbb{C}}
\def\R{\mathbb{R}}
\def\Q{\mathbb{Q}}
\def\P1{\mathbb{P}^{1}}
\def\Etproj{\overline{E}_{t}}
\def\iup{{\widetilde{\iota}}}
\def\rx{r_{x}}
\def\ry{r_{y}}

\def\mathfrakD{\mathfrak{D}}
\def\PHI{\mathcal{O}_1} 
\def\PHI{\Phi}

\author[T.~Dreyfus]{Thomas Dreyfus\textsuperscript{1}\protect\orcid{0000-0003-1459-8456}}
\address{\textsuperscript{1}Institut de Recherche Math\'ematique Avanc\'ee, U.M.R. 7501 Universit\'e de Strasbourg et C.N.R.S. 7, rue Ren\'e Descartes 67084 Strasbourg, France; \newline
\website{https://sites.google.com/site/thomasdreyfusmaths/}}

\title[Differential algebraic generating functions of walks in the quarter plane]{Differential algebraic generating functions\\ of weighted walks in the quarter plane}

\abstract{In the present paper we study the nature of the trivariate generating functions of weighted walks in the quarter plane. 
Combining the results of this paper to previous ones, we complete the proof of the following theorem. 
The series satisfies a nontrivial algebraic differential equation in one of its variables,  if and only if it satisfies a nontrivial algebraic differential equation in each of its variables.}


\begin{document}
\maketitle

\section{Introduction}

\textbf{Framework.}
Consider a  walk with small steps in the nonnegative quadrant $\mathbb{Z}_{\geq0}^{2}=\{0,1,2,\ldots\}^2$ starting from $P_0:=(0,0)$, that is a succession of points
\begin{equation*}
     P_{0}, P_1,\ldots ,P_k,
\end{equation*}   
where each $P_n$ lies in the quarter plane, where the moves (or steps) $P_{n+1}-P_{n}$ belong to $\{0,\pm 1\}^{2}$, and the probability to move in the direction $P_{n+1}-P_{n}=(i,j)$ may be interpreted as  some  weight parameter $d_{i,j}\in [0,1]$, with $\sum_{(i,j)\in\{0,\pm 1\}^{2}}d_{i,j}=1$. The model of the walk (or model for short) is the data of the $d_{i,j}$ and the step set of the  walk  is  the set of directions with nonzero  weights, that is 
$$
\mathcal{S}=\{ (i,j ) \in \{0,\pm 1\}^{2}  \ | \ d_{i,j} \neq 0\}.
$$
If $d_{0,0}=0$ and if the nonzero $d_{i,j}$ all have the same value, we say that  the model  is unweighted.
 The following picture  provides an example of a walk in the nonnegative quadrant:
\bigskip
\begin{center}
\begin{tikzpicture}[scale=.99, baseline=(current bounding box.center)]
\draw (0,0)--(10,0);
\draw (0,0)--(0,4);
\draw[->](0,0)--(1,1);
\draw[->](1,1)--(1,0);
\draw[->](1,0)--(0,1);
\draw[->](0,1)--(1,2);
\draw[->](1,2)--(2,1);
\draw[->](2,1)--(2,0);
\draw[->](2,0)--(3,1);
\draw[->](3,1)--(3,0);
\draw[->](3,0)--(4,1);
\draw[->](4,1)--(3,2);
\draw[->](3,2)--(2,3);
\draw[->](2,3)--(2,2);
\draw[->](2,2)--(3,3);
\draw[->](3,3)--(4,2);
\draw[->](4,2)--(4,1);
\draw[->](4,1)--(5,0);
\draw[->](5,0)--(6,1);
\draw[->](6,1)--(6,0);
\draw[->](6,0)--(7,1);
\draw[->](7,1)--(8,0);
\draw[->](8,0)--(9,1);
\draw[->](9,1)--(9,0);
\draw[->](9,0)--(10,1);
\draw[->](10,1)--(9,2);
\draw[->](9,2)--(8,3);
\draw[->](8,3)--(8,2);
\end{tikzpicture}
\quad\quad
$\mathcal{S}=\left\{\begin{tikzpicture}[scale=.4, baseline=(current bounding box.center)]
\draw[thick,->](0,0)--(-1,1);
\draw[thick,->](0,0)--(1,1);
\draw[thick,->](0,0)--(1,-1);
\draw[thick,->](0,0)--(0,-1);
\end{tikzpicture}\right\}$
\end{center}
\bigskip
 Such objects are very natural both in combinatorics and probability theory: they are interesting for themselves and also because they are strongly related to other discrete structures; see~\cite{BMM,DeWa-15} and references therein.\smallskip
 
The weight of the walk is defined to be the product of the weights of its component steps. For any $(i,j)\in \Z_{\geq 0}^{2}$ and any $k\in \Z_{\geq 0}$, we let $q_{i,j,k}$ be the {sum of the weights of all walks reaching}  the position $(i,j)$ from the initial position $(0,0)$ after $k$ steps. We introduce the corresponding trivariate generating function
\begin{equation}
Q(x,y;t):=\displaystyle \sum_{i,j,k\geq 0}q_{i,j,k}x^{i}y^{j}t^{k}.
\end{equation}
Being the generating function of probabilities, $Q(x,y;t)$ converges for all $(x,y,t)\in \C^{3}$ such that $\vert x\vert,\vert y\vert \leq 1$ and $\vert t\vert < 1$. Note that in several papers, as in~\cite{BMM}, it is not assumed that $\sum d_{i,j}=1$. However, after a rescaling of the $t$ variable, we may always reduce to this case.
\bigskip
 
 \textbf{Statement of the main result.}
As we will see in the sequel, this paper takes part in a long history of articles that study the algebraic and differential relations satisfied by $Q(x,y;t)$. For any choice of a variable  $\star$ among $x,y,t$, we say that 
 $Q(x,y;t)$ is $\partial_{\star}$-algebraic if 
 there exists $n\in \Z_{\geq 0}$, such that there exists a nonzero multivariate polynomial $P_{\star}\in  \C(x,y,t)[X_{0},\dots,X_{n}]$, such that 
\begin{equation}
0=P_{\star}(Q(x,y;t),\dots,\partial_{\star}^{n} Q(x,y;t)).
\end{equation}
   We stress that in the above definition, it is equivalent to  require 
   ${0\neq P_{\star} \in  \Q[X_{0},\dots,X_{n}]}$; see Remark~\ref{rem2}. 
   Otherwise, we say that the series $Q(x,y;t)$ is  $\partial_{\star}$-differentially transcendental.\smallskip
   
 Since the three  variables $x,y$ and $t$ play a different role, we might expect the series to be of different nature with respect to the three derivatives. The main result of this paper, quite unexpected at first sight, shows that it is not the case. More precisely, using results of this paper and combining them to partial cases already known (see the discussion in the sequel), we complete the proof of the following main theorem. 

\begin{theorem}\label{thm1}
The following facts are equivalent:
\begin{itemize}
\item The series $Q(x,y;t)$ is $\partial_{x}$-algebraic;
\item The series $Q(x,y;t)$ is $\partial_{y}$-algebraic; 
\item The series $Q(x,y;t)$ is $\partial_{t}$-algebraic. 
\end{itemize}
\end{theorem}
Note that an algorithm is given in~\cite[Section 5]{hardouin2020differentially} to decide whether the generating function is differentially algebraic in the $x$ variable or not, but this does not provide the differential equation when it exists.
\bigskip

\textbf{State of the art.}
More generally, the question of studying whether $Q(x,y;t)$ satisfies  algebraic (resp.~linear differential, resp.~algebraic differential) equations attracted the attention of many authors in the last decade. In the unweighted case, the problem was first addressed in the seminal paper~\cite{BMM} and solved using several methods, such as combinatorics, computer algebra, complex analysis, and more recently, difference Galois theory; see~\cite{mishna2009two,BostanKauersTheCompleteGenerating,kurkova2012functions,
melczer2014singularity,
DHRS,dreyfus2019length,dreyfus2020walks,bernardi2017counting}. We refer to the introduction of~\cite{dreyfus2020nature} for a history of the cited results, from which it follows that Theorem~\ref{thm1} is valid for the  unweighted models. 

\pagebreak

The main difficulty in generalizing those results to weighted models is that, contrary to the unweighted framework, there are  infinitely many  weighted models.
 However, certain unweighted results  are still valid in the weighted cases, while some others are proved by a case-by-case argument, and therefore cannot be generalized straightforwardly.  So beyond the generalization, we believe that replacing case-by-case proofs by systematic arguments has its own interest since it shows that the unweighted version of Theorem~\ref{thm1} has not appeared by accident in a finite number of cases, and illustrates a general phenomenon.
\smallskip 

In many situations, the equivalence between the $\partial_{x}$-algebraicity and the $\partial_{y}$-algebraicity can be straightforwardly deduced in this weighted context from the proof of~\cite[Proposition~3.10]{DHRS}. In~\cite[Theorem 2]{dreyfus2019length} it was proved that the $\partial_{t}$-algebraicity implies the $\partial_{x}$-algebraicity. So it remains to show the converse. In~\cite{bernardi2017counting}, the authors show that all $\partial_x$-differentially algebraic unweighted models have a decoupling function. They use this property to prove the $\partial_t$-algebraicity in that case. In~\cite{DHRS}, using difference Galois theory, the authors show that such unweighted models admit a telescoping relation. We refer to~\cite{hardouin2020differentially} for precise definitions of the two latter notions. In 
\cite{dreyfus2017differential}, it is proved that the  $\partial_{x}$-algebraic weighted models also have a telescoping relation.  Finally in~\cite{hardouin2020differentially} the equivalence between the existence of a telescoping relation and the existence of decoupling functions is shown. This implies that a $\partial_{x}$-algebraic series admits a certain decomposition into elliptic functions.  \smallskip 

The main difficulty is that the existence of such decompositions is proved for fixed values of $t$, so nothing is known about the dependence in $t$ of the coefficients. 
For instance,   the function $x\Gamma(t)$, seen as a function of $x$ is simple for all fixed value of $t$ (it is rational!)   but it is differentially transcendental with respect to $t$, due to  H\"older's result.  We then have to make a  careful study of the $t$-dependence of such elliptic relations, and use  some results of $\partial_t$-algebraicity of the Weierstrass function in~\cite{bernardi2017counting}.  Finally, we are able to show that 
the $\partial_{x}$-algebraicity implies the $\partial_{t}$-algebraicity. The following diagram summarizes the various contributions  toward the proof of Theorem~\ref{thm1}. 
\begin{equation*}
\xymatrix{
    \partial_y\hbox{-algebraicity} \ar@{<=>}[r]_{\hbox{\cite{DHRS}}} &\partial_x\hbox{-algebraicity }\ar@{<=}[r]_{\hbox{\cite{dreyfus2019length}}} &\partial_t\hbox{-algebraicity}  \\
     & \hbox{Telescoping relation} \ar@{<=}[u]_{\hbox{\cite{dreyfus2017differential}}}\ar@{<=>}[r]^{\hbox{\cite{hardouin2020differentially}}}&\hbox{Decoupling function}\ar@{=>}[u]_{\hbox{This paper +\cite{bernardi2017counting}}}
  }  
\end{equation*}
\smallskip
 
\textbf{Structure of the paper.}
The paper is organized as follows. In Section~\ref{sec2} we provide some reminders of objects appearing in the study of models of walks in the quarter plane. More precisely, we will study well-known properties of the kernel curve and explain  how the generating function may be continued. We will also explain why Theorem~\ref{thm1} is correct in some degenerate cases that we may withdraw.
In Section~\ref{sec3} we prove technical results on differential algebraicity.  Some intermediate results stay valid in the framework of algebraic functions and/or solution of linear differential equations, but to simplify the exposition,  we chose to present this section in a unified framework, making some intermediate results suboptimal.  Finally Section~\ref{sec4} is devoted to the proof of Theorem~\ref{thm1}. We split our study in two cases depending on whether the so-called group of the walk is finite or not.
\bigskip
\pagebreak

\section{Kernel of the walk}\label{sec2}

\subsection{Functional equation}\label{sec21}
The kernel of the walk is the polynomial defined by 
$$K(x,y;t):=xy (1-t S(x,y)),$$ 
where $S(x,y)$ denotes the jump polynomial
\begin{eqnarray*}
S(x,y) &=&\sum_{(i,j)\in \{0,\pm 1\}^{2}} d_{i,j}x^i y^j \\
&=& A_{-1}(x) \frac{1}{y} +A_{0}(x)+ A_{1}(x) y \\ 
&=&  B_{-1}(y) \frac{1}{x} +B_{0}(y)+ B_{1}(y) x,
\end{eqnarray*}
with $A_{i}(x) \in x^{-1}\R[x]$, $B_{i}(y) \in y^{-1}\R[y]$ (we recall that we consider weights $d_{i,j} \in [0,1]$).
The kernel plays an important role in the so-called kernel method and the techniques we are going to apply will vary depending on its algebraic properties, that have been studied in~\cite{FIM} (when $t=1$), and in~\cite{dreyfus2017differential,DHRS20} (when $t\in (0,1)$). 
The starting point is the following fundamental functional equation.


\begin{lemma}\label{lem1}
The generating function $Q(x,y;t)$  satisfies the functional equation
\begin{equation}
     K(x,y;t)Q(x,y;t)=xy+
     K(x,0;t)Q(x,0;t) +K(0,y;t)Q(0,y;t)-K(0,0;t) Q(0,0;t).
\end{equation}
\end{lemma}
\begin{proof}
As a walk is either empty, or a smaller walk to which one added a step (removing the cases leaving the quarter-plane), one has the following combinatorial functional equation
\begin{equation}
\resizebox{\linewidth}{!}{
$Q(x,y;t)=1+ t S(x,y) Q(x,y;t) 
- t  \frac{B_{-1}(y)}{x}   Q(0,y;t) 
- t \frac{A_{-1}(x)}{y}   Q(x,0;t) 
+ t \frac{d_{-1,-1}}{xy} Q(0,0;t),
$}
\end{equation}
\vspace{-8mm}

\noindent where the last summand removes the corresponding double counting. 
Multiplying by~$xy$, we get Lemma~\ref{lem1}.
\end{proof}

\subsection{Degenerate cases.}

Like in~\cite{DHRS20}, we will discard the following degenerate cases.
 
\begin{definition}[Degenerate model]
Let us fix $t\in (0,1)$.  A  model of walk is called degenerate if one of the following holds:
\begin{itemize}
\item $K(x,y;t)$ factors in non-constant polynomials in  $\C[x,y]$, 
\item $K(x,y;t)$ has $x$-degree (or $y$-degree) less than or equal to $1$.
\end{itemize}
 \end{definition}
 
The notion of degeneracy thus seems to depend upon the parameter $t$.    However,  we will see in Proposition~\ref{prop:singcases} below that the model is degenerate for a value of $t\in (0,1)$ if and only if it is degenerate for all values of  $t\in (0,1)$.  So, to lighten the terminology, we prefer not to stress this  $t$-dependence and we say ``degenerate'' rather than ``$t$-degenerate''. 

In what follows we will sometimes represent a family of models of walks with arrows.  For instance, the family of models  represented by 
$\begin{tikzpicture}[scale=0.357, baseline=(current bounding box.center)]
\draw[thick,->](0,0)--(-1,-1);
\draw[thick,->](0,0)--(1,-1);
\draw[thick,->](0,0)--(1,1);
\draw[thick,->](0,0)--(0,1);
\end{tikzpicture} \hbox{ \text{\  or, equivalently, \   } } \left\{\begin{tikzpicture}[scale=.4]
\draw[thick,->](0,0)--(1,1);
\end{tikzpicture},\ \begin{tikzpicture}[scale=.4]
\draw[thick,->](0,0)--(1,-1);
\end{tikzpicture},\ \begin{tikzpicture}[scale=.4]
\draw[thick,->](0,0)--(-1,-1);
\end{tikzpicture}, \ 
\begin{tikzpicture}[scale=.4)]
\draw[thick,->](0,0)--(0,1);
\end{tikzpicture} \ 
\right\}$
corresponds to models with $d_{1,0}=d_{0,-1}=d_{-1,1}=d_{-1,0}=0$ and nothing is assumed on the other $d_{i,j}$. We stress the fact that the other $d_{i,j}$ (the weight of the arrows above) are allowed to be $0$. In the following results, 
 the behavior of the kernel curve \textit{never} depends on $d_{0,0}$. This is the reason why, to reduce the amount of notations, we have decided  not to associate an arrow to $d_{0,0}$.\smallskip 
 
The following proposition has been proved in~\cite[Lemma 2.3.2]{FIM} for $t=1$, in~\cite[Proposition 1.2]{DHRS20} for $t$ is transcendental over $\Q(d_{i,j})$, and in~\cite[Proposition 1.3]{dreyfus2017differential} for the other values of $t$ in $(0,1)$.

 \begin{proposition}\label{prop:singcases}
Let us fix $t\in (0,1)$. A  model of walk is  degenerate if and only if at least one of the following holds: 
\begin{enumerate}[label=\textnormal{({\alph*})},resume]
\item \label{case1}There exists $i\in \{- 1,1\}$ such that $d_{i,-1}=d_{i,0}=d_{i,1}=0$. This corresponds to the following  families of   models
$$\begin{tikzpicture}[scale=0.6, baseline=(current bounding box.center)]
\draw[thick,->](0,0)--(0,-1);
\draw[thick,->](0,0)--(1,-1);
\draw[thick,->](0,0)--(1,0);
\draw[thick,->](0,0)--(1,1);
\draw[thick,->](0,0)--(0,1);
\end{tikzpicture} \ \ ,\quad 
\begin{tikzpicture}[scale=0.6, baseline=(current bounding box.center)]
\draw[thick,->](0,0)--(0,-1);
\draw[thick,->](0,0)--(-1,-1);
\draw[thick,->](0,0)--(-1,0);
\draw[thick,->](0,0)--(-1,1);
\draw[thick,->](0,0)--(0,1);
\end{tikzpicture}$$ 

\item \label{case2} There exists $j\in \{-1, 1\}$ such that $d_{-1,j}=d_{0,j}=d_{1,j}=0$. This corresponds to the following   families of   models
$$\begin{tikzpicture}[scale=.6, baseline=(current bounding box.center)]
\draw[thick,->](0,0)--(-1,0);
\draw[thick,->](0,0)--(-1,-1);
\draw[thick,->](0,0)--(0,-1);
\draw[thick,->](0,0)--(1,-1);
\draw[thick,->](0,0)--(1,0);
\end{tikzpicture} \ \ ,
 \quad
\begin{tikzpicture}[scale=.6, baseline=(current bounding box.center)]
\draw[thick,->](0,0)--(-1,0);
\draw[thick,->](0,0)--(-1,1);
\draw[thick,->](0,0)--(0,1);
\draw[thick,->](0,0)--(1,1);
\draw[thick,->](0,0)--(1,0);
\end{tikzpicture}$$

\item \label{case3} All the  weights are $0$ except maybe  $\{d_{-1,-1},d_{0,0},d_{1,1}\}$ or  $\{d_{-1,1},d_{0,0},d_{1,-1}\}$. This corresponds to the following  families of  models
$$
\left\{\begin{tikzpicture}[scale=.6, baseline=(current bounding box.center)]
\draw[thick,->](0,0)--(1,1);
\end{tikzpicture}, \begin{tikzpicture}[scale=.6, baseline=(current bounding box.center)]
\draw[thick,->](0,0)--(-1,-1);
\end{tikzpicture} \right\} \ ,
\quad
\left\{\begin{tikzpicture}[scale=.6, baseline=(current bounding box.center)]
\draw[thick,->](0,0)--(-1,1);
\end{tikzpicture}, \begin{tikzpicture}[scale=.6, baseline=(current bounding box.center)]
\draw[thick,->](0,0)--(1,-1);
\end{tikzpicture} \right\} 
$$
\end{enumerate}
\end{proposition} 

In virtue of the following lemma, Theorem~\ref{thm1} is valid for the degenerate models of walks. Therefore we will focus on models that are not degenerate.
 
\begin{lemma}\label{lem2}
Assume that the model of walk is  degenerate. Then $Q(x,y;t)$ is algebraic over $\C(x,y,t)$ (and thus is differentially algebraic in its three variables).
\end{lemma}

\begin{proof}
We use Proposition~\ref{prop:singcases}.
Consider the cases~\ref{case1},~\ref{case2}, and first configuration of the case~\ref{case3}.
In the unweighted case  it is proved in~\cite[Section 1.2]{BMM} that $Q(x,y;t)$ is algebraic over $\C(x,y,t)$. The proof is the same in the weighted context but, to be self-contained, let us sketch the proof here. In the first configuration of case~\ref{case1}  the generating function is the same as the  corresponding  generating function of a model in the upper half-plane $\Z\times \N$. The latter is classically known to be algebraic over $\C(x,y,t)$,  see for instance~\cite[Proposition~2]{bousquet2003walks}.  In the second configuration of case~\ref{case1}, we have a unidimensional walk on the $y$-axis and such series is known to be rational, and therefore algebraic over $\C(x,y,t)$. The case~\ref{case2} is similar. In the first configuration of case~\ref{case3}, we are considering a unidimensional walk on the half-line $\{(x,x),x\in \N\}$, and the generating function is algebraic. 
Since in all these cases,  $Q(x,y;t)$ is algebraic over $\C(x,y,t)$, it is differentially algebraic in its three variables.
In the last configuration of case~\ref{case3}, all the  weights are $0$ except maybe  $\{d_{-1,1},d_{0,0},d_{1,-1}\}$, so the walk cannot leave $(0,0)$ and we have 
\begin{equation}\label{Q00}
 Q(x,y;t)= \sum_{k=0}^{\infty}d_{0,0}^{k}t^{k}=\frac{1}{1-d_{0,0}t}.
\end{equation}
Therefore the result holds in that case too.
\end{proof}
\pagebreak 

\subsection{Genus of the walk}\label{Genus}

The  kernel curve $E_{t}$ is the complex affine algebraic curve defined~by
\begin{equation*}
 E_t = \{(x,y) \in \C \times \C \ \vert \ K(x,y;t) = 0\}.
 \end{equation*} 
We shall now consider a compactification of this curve. We let $\P1(\C)$ be the complex projective line, that is the quotient of $(\C \times \C) \setminus \{(0,0)\}$ by the equivalence relation $\sim$ defined by 
\begin{equation}
(x_{0},x_{1}) \sim (x_{0}',x_{1}') \Leftrightarrow \exists \lambda \in \C^{*},  (x_{0}',x_{1}') = \lambda (x_{0},x_{1}). 
\end{equation}
The equivalence class of $(x_{0},x_{1}) \in (\C \times \C) \setminus \{(0,0)\}$ is denoted by $[x_{0}:x_{1}] \in \P1(\C)$. The map 
$x \mapsto  [x:1]$  embeds $\C$ inside $\P1(\C)$. The latter map is not surjective: its image is $\P1(\C) \setminus \{[1:0]\}$; the missing point $[1:0]$  is usually denoted by $\infty$. 
  Now, we embed $E_{t}$  inside $\P1(\C) \times \P1(\C)$ via  ${(x,y) \mapsto ([x:1],[y:1])}$. The  kernel curve $\Etproj$ is the closure of this embedding of $E_{t}$.  In other words, the  kernel curve $\Etproj $ is the algebraic curve defined by
\begin{equation}
\Etproj = \{([x_{0}:x_{1}],[y_{0}:y_{1}]) \in \P1(\C) \times \P1(\C) \ \vert \ \overline{K}(x_0,x_1,y_0,y_1;t) = 0\}
\end{equation}
where $\overline{K}(x_0,x_1,y_0,y_1;t)$ is the following degree two homogeneous polynomial in the four variables $x_0,x_1,y_0,y_1$
\begin{equation*}
\overline{K}(x_0,x_1,y_0,y_1;t)=x_1^2y_1^2K\left(\frac{x_0}{x_1},\frac{y_0}{y_1};t\right)= x_0x_1y_0y_1 -t \sum_{i,j=0}^2 d_{i-1,j-1} x_0^{i} x_1^{2-i}y_0^j y_1^{2-j}. 
 \end{equation*}

Although it may seem more natural to take the closure of $\Etproj$ in $\mathbb{P}^2(\C)$, the above definition allows us to avoid unnecessary singularities.\smallskip

The following proposition has been proved in~\cite[Proposition 2.1 and Corollary~2.6]{DHRS20}, when $t$ is transcendental over $\Q(d_{i,j})$ and has been extended for a general $0<t<1$ in~\cite[Proposition 1.9]{dreyfus2017differential}.

\begin{proposition}\label{prop2}
Let us fix $t\in (0,1)$ and assume that the model of the walk is not degenerate. The following facts are equivalent:
\begin{enumerate}
\item $\Etproj$ is an elliptic curve;
\item  The set of authorized directions $\mathcal{S}$ is not included in any half-space with boundary passing through the origin.
\end{enumerate}
\end{proposition}

Let us now discuss the case where for $t\in (0,1)$ fixed, the model  is not degenerate and $\Etproj$ is not an elliptic curve. 
By Proposition~\ref{prop:singcases} and Proposition~\ref{prop2}, this corresponds to nondegenerate models that belong to one of the four families in Figure~\ref{4models}.
\begin{figure}[hb]
\begin{equation}
\begin{tikzpicture}[scale=.6,baseline=(current bounding box.center)]
\draw[thick,->](0,0)--(-1,1);
\draw[thick,->](0,0)--(0,1);
\draw[thick,->](0,0)--(1,1);
\draw[thick,->](0,0)--(1,0);
\draw[thick,->](0,0)--(1,-1);
\end{tikzpicture}\quad 
\begin{tikzpicture}[scale=.6,baseline=(current bounding box.center)]
\draw[thick,->](0,0)--(-1,1);
\draw[thick,->](0,0)--(1,1);
\draw[thick,->](0,0)--(-1,0);
\draw[thick,->](0,0)--(0,1);
\draw[thick,->](0,0)--(-1,-1);
\end{tikzpicture}\quad
\begin{tikzpicture}[scale=.6,baseline=(current bounding box.center)]
\draw[thick,->](0,0)--(-1,1);
\draw[thick,->](0,0)--(-1,0);
\draw[thick,->](0,0)--(-1,-1);
\draw[thick,->](0,0)--(0,-1);
\draw[thick,->](0,0)--(1,-1);
\end{tikzpicture}\quad
\begin{tikzpicture}[scale=.6,baseline=(current bounding box.center)]
\draw[thick,->](0,0)--(1,1);
\draw[thick,->](0,0)--(1,0);
\draw[thick,->](0,0)--(-1,-1);
\draw[thick,->](0,0)--(0,-1);
\draw[thick,->](0,0)--(1,-1);
\end{tikzpicture} \qquad
\end{equation}
\caption{Our four nondegenerate models}\label{4models}
\end{figure}

Note that although the third configuration in~Figure~\ref{4models} is called nondegenerate, it leads to walks that never escape from $(0,0)$ and thus their generating function is trivial. 

The following lemma yields that Theorem~\ref{thm1} is valid for the families of models in Figure~\ref{4models}.
\pagebreak

\begin{lemma}\label{lem3}
The following holds.
\begin{enumerate}[label=\textnormal{(\alph*)}]
\item Assume that the model of the walk is not degenerate and  belongs to the first family in Figure~\ref{4models}. Then $Q(x,y;t)$ is differentially transcendental in its three variables. 
\item Assume that the model of the walk belongs to the second, third or the fourth family in Figure~\ref{4models}. Then $Q(x,y;t)$ is algebraic over $\C(x,y,t)$, and thus is differentially algebraic in its three variables.
\end{enumerate}
\end{lemma}
\begin{proof}
\begin{enumerate}[label=\textnormal{(\alph*)}]
\item This is~\cite[Corollary~2.2]{dreyfus2019length}; see also~\cite[Theorem 3.1]{dreyfus2020walks}.
\item Consider the second family.  We have 
$Q(x,0;t)=Q(0,0;t)$ and $K(x,0;t)=K(0,0;t)$. Then by Lemma~\ref{lem1}, 
\begin{equation}\label{eq8}
 K(x,y;t)Q(x,y;t)=K(0,y;t)Q(0,y;t)+xy.
\end{equation}
\noindent Let us see that with the same arguments as for the walks in the half-plane,  we deduce that $Q(x,y;t)$ is  algebraic over $\C(x,y,t)$.  The idea is to   locally write   $K(\phi(y;t),y;t)=0$.   Evaluating at $(\phi(y;t),y;t)$ we then have for convenient $y$ and $t$,  $0=K(0,y;t)Q(0,y;t)+\phi(y;t)y$, proving that  $Q(0,y;t)$ is algebraic over $\C(x,y,t)$.  The functional equation~\eqref{eq8} allows then to conclude that $Q(x,y;t)$ is algebraic over $\C(x,y,t)$.
As in the proof of  Lemma~\ref{lem2}, we may deduce that  $Q(x,y;t)$ is differentially algebraic in its three variables.
 The reasoning for the fourth family is similar. For the third family, the walk has to stay at $(0,0)$ and we have $$Q(x,y;t)=\displaystyle \sum_{k=0}^{\infty}d_{0,0}^{k}t^{k}=\frac{1}{1-d_{0,0}t}.$$
Therefore the result holds in that case too.\qedhere
\end{enumerate}
\end{proof}

\subsection{Group of the walk}
From now on, we may focus on the case where $\Etproj$ is an elliptic curve.  Recall that we have seen in Proposition~\ref{prop:singcases}, that  $K(x,y;t)$ has degree  two in $x$ and~$y$, and nonzero coefficient of degree $0$ in $x$ and $y$. Hence, $A_{1}(x),A_{-1}(x),B_{1}(y),B_{-1}(y)$ are not identically zero. \smallskip

Following~\cite[Section 3]{BMM},~\cite[Section 3]{KauersYatchak}  or~\cite{FIM}, and using the notations introduced in Section~\ref{Genus}, we consider the rational involutions 
 given by 
\begin{equation*}
i_1([x_0:x_1],[y_0:y_1]) =\left(\frac{x_{0}}{x_{1}}, \frac{A_{1}(\frac{x_{0}}{x_{1}}) }{A_{-1}(\frac{x_{0}}{x_{1}})\frac{y_{0}}{y_{1}}}\right)
\text{\ \ and \ \ }  i_2([x_0:x_1],[y_0:y_1])=\left(\frac{B_{-1}(\frac{y_{0}}{y_{1}})}{B_{1}(\frac{y_{0}}{y_{1}})\frac{x_{0}}{x_{1}}},\frac{y_{0}}{y_{1}}\right).
\end{equation*}
Note that we have $i_1([x_0/x_1:1],[y_0/y_1:1])=i_1([x_0:x_1],[y_0:y_1])$ and the same holds for~$i_2$.
Note also that $i_{1}, i_{2}$  are  a priori not defined when the denominators vanish but we will see in the sequel that we may overcome this problem when we will restrict to $\Etproj$. \smallskip

For a fixed value of $x$, there are at most  two possible values of $y$ such that $(x,y)\in\Etproj$. The involution $i_{1}$ corresponds to interchanging these values. A similar interpretation can be given for~$i_2$. Therefore  the kernel curve $\Etproj$ is left invariant by the natural action of~$i_{1}, i_{2}$. \pagebreak

\begin{figure}
\begin{center}
\includegraphics{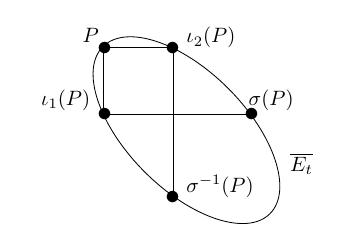}
\vspace{-5mm}
\end{center}
\caption{The maps $i_{1}$ and $i_{2}$ restricted to the kernel curve $\Etproj$ are denoted by $\iota_1$ and $\iota_2$, respectively.}\label{figiota}
\end{figure}

 We denote by $\iota_{1}, \iota_{2}$ the restriction of $i_{1},i_{2}$ to $\Etproj$; see Figure~\ref{figiota}. Again, these functions are  a priori not defined where the denominators vanish. However, by~\cite[Proposition~4.1]{DHRS20}, this is only an ``apparent problem''. To be precise, the authors  proved this for $t$ transcendental over $\Q(d_{i,j})$ but the proof is still valid  when $\Etproj$ is an elliptic curve. 
 We then obtain that
 $\iota_{1}$ and $\iota_{2}$ can be extended to morphisms of $\Etproj$. We recall that a rational map $f$
 from $\Etproj$ to~$\Etproj$ is a morphism if it is regular at any $P \in \Etproj$,  i.e.~if $f$ can be represented in suitable affine charts containing $P$ and  $f(P)$ by a rational function with nonvanishing denominator at $P$. \smallskip
 
Let us finally define $\sigma=\iota_2 \circ \iota_1$.
Note that such a map is called a QRT-map and has been widely studied; see~\cite{DuistQRT}. 

\begin{definition}[Group of the walk]\label{group}
We call $G$  the group generated by $\iota_{1}$ and $\iota_{2}$ and we call $G_{t}$ the specialization of this group for a fixed value of $0<t<1$. 
\end{definition}

 In the unweighted case,  the algebraic nature of the generating  series depends on whether~$\sigma$ has finite or infinite order. More precisely, $G$ is finite if and only if the generating function is holonomic, i.e.~satisfies a nontrivial linear differential equation with coefficients in $\C(x,y,t)$ in each of its three variables.  On the other hand, when $G$ is infinite, $G_{t}$ can be either finite or infinite; see~\cite{fayolleRaschel} for concrete examples. However, in that situation, the set of values of $t$ such that $G_t$ is finite is countable,  see~\cite[Proposition~2.6]{DHRS}.

\subsection{Uniformization of the curve}\label{sec15}
In this section, we consider the uniformization problem in the genus one context, that has been solved in~\cite{FIM} for the case $t=1$, and 
~\cite{dreyfus2017differential} for the case $0<t<1$. 
Let us consider a  nondegenerate model of walk and assume that for all $t\in (0,1)$, $\Etproj$ is an elliptic curve. 
 By Proposition~\ref{prop2}, this corresponds to the situation where the step  set is not included in any half-plane whose boundary passes through $(0, 0)$.
By~\cite[Proposition 2.1]{dreyfus2017differential}, the elliptic curve $\Etproj$ is biholomorphic to $\C/(\omega_1 (t)\Z  + \omega_2 (t)\Z)$ for some lattice $\omega_1 (t)\Z+ \omega_2 (t)\Z$ of $\C$  via some $(\omega_1 (t)\Z + \omega_2 (t)\Z)$-periodic holomorphic map $\Lambda$ 
\begin{equation}
\begin{array}{llll}
\Lambda :& \C& \rightarrow &\overline{E}_t\\
 &\Lambda(\omega) &:=& (x(\omega;t), y(\omega;t)),
\end{array}\label{Lambda}
\end{equation}
where $x, y$ are rational functions of $\wp$ and its derivative $\partial_{\omega}\wp$, and $\wp$ is the Weierstrass function associated with the lattice $\omega_1 (t)\Z + \omega_2 (t)\Z$:
\begin{equation*}
     \wp(\omega ;t)=\frac{1}{\omega^{2}}+ \sum_{(\ell_{1},\ell_{2}) \in \Z^{2}\setminus \{(0,0)\}} \left(\frac{1}{(\omega +\ell_{1}\omega_{1}(t)+\ell_{2}\omega_{2}(t))^{2}} -\frac{1}{(\ell_{1}\omega_{1}(t)+\ell_{2}\omega_{2}(t))^{2}}\right).
\end{equation*}     
Then, the field of meromorphic functions on $\Etproj$ is isomorphic to the field of meromorphic functions on $\C/(\omega_1 (t)\Z+ \omega_2 (t)\Z)$, that is itself isomorphic to the field of meromorphic functions on $\C$ that are $(\omega_{1} (t),\omega_{2} (t))$-periodic. For $t\in (0,1)$ fixed, the latter field  is equal to $\C(\wp, \partial_{\omega}\wp )$; see~\cite[Chapter~9, Theorem 1.8]{stein2010complex}.
 
The maps $\iota_{1}$, $\iota_{2}$, and $\sigma$ may be analytically lifted to the $\omega$-plane $\C$ via the map $\Lambda^{-1}$. We denote these lifts by  $\iup_{1}$, $\iup_{2}$, and $\widetilde{\sigma}$ respectively. So we have the commutative diagrams 
\begin{equation*}
\xymatrix{
    \Etproj  \ar@{->}[r]^{\iota_k} & \Etproj  \\
    \C \ar@{->}[u]^\Lambda \ar@{->}[r]_{\iup_k} & \C \ar@{->}[u]_\Lambda 
  }
  \qquad\qquad\qquad
  \xymatrix{
    \Etproj  \ar@{->}[r]^{\sigma} & \Etproj  \\
\C \ar@{->}[u]^\Lambda \ar@{->}[r]_{\widetilde{\sigma}} & \C \ar@{->}[u]_\Lambda 
  }  
\end{equation*}
 For any $[x_0:x_1]$ in $\P1(\C)$, we denote by $\Delta_{1}([x_0:x_1];t)$ the discriminant of the degree two homogeneous polynomial given by $y \mapsto \overline{K}(x_0,x_1,y;t)$.
 Let us write 
$$\Delta_{1}([x_0:x_1];t)= \displaystyle\sum_{i=0}^{4}\alpha_{i}(t)x_{0}^{i}x_{1}^{4-i}.$$
 By~\cite[Theorem 1.11]{dreyfus2017differential}, the discriminant $\Delta_{1}([x_0:x_1];t)$ admits four distinct continuous real roots $a_1(t),\dots,a_{4}(t)$. They are numbered such that the cycle of $\mathbb{P}^{1}(\R)$ starting from $-1$ to~$\infty$ and from $-\infty$ to $-1$ crosses the $a_i$ in the order $a_1 (t),\dots, a_{4}(t)$. 
Furthermore, $[1:0]$ is a root  if and only if $\alpha_{4}(t)=0$. In~\cite[Section 1.4]{dreyfus2017differential}, we see that $\alpha_{4}(t)=t^{2}(d_{1,0}^{2}-4d_{1,-1}d_{1,1})$. It follows that $[1:0]$ is a root of $\Delta_{1}([x_0:x_1];t)$ for one value of $t\in (0,1)$, if and only if $[1:0]$ is a root of $\Delta_{1}([x_0:x_1];t)$ for all  $t\in (0,1)$. \smallskip

 Similarly, we denote by  $b_1(t),\dots,b_{4}(t)$ the continuous real roots of the discriminant $x \mapsto \overline{K}(x,y_0,y_1;t)$, numbered in the same way, and we write ${\Delta_{2}([y_0:y_1];t)= \displaystyle\sum_{i=0}^{4}\beta_{i}(t)y_{0}^{i}y_{1}^{4-i}}$.

The following formulas have been proved
\begin{itemize}
\item in~\cite[Section 3.3]{FIM} when $t=1$,
\item in~\cite{RaschelJEMS} in the unweighted case,
\item in~\cite[Proposition 2.1 and (2.16)]{dreyfus2017differential}, in the weighted case.
\end{itemize}

\begin{proposition}[\!\!\cite{dreyfus2017differential}, Proposition 2.1, Lemma 2.6, and (2.16)]
\label{prop:uniformization}
For $i=1,2$, let us set $D_i(\star;t):=\Delta_{i}([\star:1];t)$. An explicit expression of the periods is given by the elliptic integrals
\begin{equation*}
     \omega_{1} (t)=\mathbf{i}\int_{a_{3}(t)}^{a_{4}(t)} \frac{dx}{\sqrt{\vert D_1(x;t)\vert}}\in \mathbf{i}\R_{>0}\quad\text{and}\quad
\omega_{2}(t)=\int_{a_{4}(t)}^{a_{1}(t)} \frac{dx}{\sqrt{D_1(x;t)}}\in \R_{>0}.
\end{equation*}
An explicit expression of the holomorphic  map $\Lambda(\omega;t)=(x(\omega;t),y(\omega;t))$ is given by 
\begin{itemize}
\item
If $a_{4}(t)\neq [1\!:\!0]$, then $x(\omega;t)=\left[a_{4}(t)+\frac{D'_{1}(a_{4}(t);t)}{\wp (\omega;t)-\frac{1}{6}D''_{1}(a_{4}(t);t)}:1\right]$;
\vspace{0.2cm}
\item  If $a_{4}(t)=[1\!:\!0]$, then $x(\omega;t)=
\left[\wp(\omega;t)-\alpha_{2}(t)/3:\alpha_{3}(t)\right]$;\vspace{0.2cm}
\item If $b_{4}(t)\neq [1\!:\!0]$, then
$y(\omega;t)=\left[b_{4}(t)+\frac{D'_{2}(b_{4}(t);t)}{\wp (\omega-\omega_3 (t)/2;t)-\frac{1}{6}D''_{2}(b_{4}(t);t)}:1\right]$;\vspace{0.2cm}
\item  If $b_{4}(t)=[1\!:\!0]$, then $y(\omega;t)=\left[\wp(\omega-\omega_3 (t)/2;t)-\beta_{2}(t)/3:\beta_{3}(t)\right]$.
\end{itemize}
An explicit expression of the involutions is given by
$$\iup_{1}(\omega)=-\omega, \quad\iup_{2}(\omega)=-\omega+\omega_{3}\quad \text{and} \quad\widetilde{\sigma}(\omega)=\omega+\omega_{3},$$
where 
\begin{equation}\label{omega3}
\omega_{3}(t)=\int_{a_{4}(t)}^{X_{\pm}(b_{4}(t);t)} \frac{dx}{\sqrt{D_1(x;t)}}\in (0,\omega_{2}(t)),
\end{equation}
and $X_{\pm}(y;t)$ are the two roots of $\overline{K}(X_{\pm}(y;t),y;t)=0$.
\end{proposition}

\subsection{Meromorphic continuation of the generating function}\label{sec16}

Let us summarize here the results of~\cite[Section 2.3]{dreyfus2017differential}. Let us fix $t\in (0,1)$.
The generating function $Q(x,y;t)$ converges for $|x|,|y|<1$. The projection of this set  inside $\mathbb{P}^{1}(\C)\times \mathbb{P}^{1}(\C)$ has a nonempty intersection with the kernel curve $\Etproj$. In virtue of~Lemma~\ref{lem1}, we then find for $|x|,|y|<1$ and $(x,y)\in \Etproj$,
$$0=K(x,0;t)Q(x,0;t)
+K(0,y;t)Q(0,y;t)
-K(0,0;t) Q(0,0;t)+xy.$$
To shorten several expressions hereafter, 
it is convenient to rewrite this equation introducing new auxiliary series $F_{1}$ and $F_{2}$:
\begin{equation}\label{eq9}
0=F_{1}(x;t) +F_{2}(y;t)-K(0,0;t) Q(0,0;t)+xy.
\end{equation}
Since the  series $F_{1}(x;t)$ and  $F_{2}(y;t)$ converge for $|x|$ and $|y|<1$ respectively,
we then continue $F_{1}(x;t)$  for $(x,y)\in \Etproj$ and $|y|<1$ with the formula: 
$$F_{1}(x;t) =-F_{2}(y;t)+K(0,0;t) Q(0,0;t)-xy.$$
We continue $F_{2}(y;t)$  for $(x,y)\in \Etproj$ and $|x|<1$ similarly. There exists a connected set $\mathcal{O}\subset \C$ such that 
\begin{itemize}
\item $\Lambda (\mathcal{O})=\{(x,y)\in \Etproj \hbox{ such that }  |x|<1 \hbox{ or } |y|<1 \}$;
\item $\widetilde{\sigma}^{-1}(\mathcal{O})\cap \mathcal{O}\neq \varnothing$;
\item $\displaystyle \bigcup_{\ell\in \Z}\widetilde{\sigma}^{\ell}(\mathcal{O})=\C$.
\end{itemize}
There also exist meromorphic functions on $\mathcal{O}$, $\rx(\omega;t)$ and $\ry(\omega;t)$, such that $\rx(\omega;t)=F_{1}(x(\omega;t);t)$  and $\ry(\omega;t)=F_{2}(y(\omega;t);t)$. 

\begin{lemma}[Inclusion of poles]\label{lem4}
The set of poles of $\rx(\omega;t)$ inside $\mathcal{O}$ are contained in the set of  poles of $x(\omega;t)$  with $|y(\omega;t)|<1$. The set of poles of $\ry(\omega;t)$ inside $\mathcal{O}$ are contained in the set of poles of $y(\omega;t)$  with $|x(\omega;t)|<1$. 
\end{lemma}

\begin{proof}
Let us use~\eqref{eq9}.  On  $\mathcal{O}$,  we have 
$$0=\rx(\omega;t) +\rx(\omega;t)-K(0,0;t) Q(0,0;t)+x(\omega;t)y(\omega;t). $$
Let us focus on $\rx(\omega;t)$,  the proof for  $\ry(\omega;t)$ is similar.  
Recall that $F_{1}(x;t)$ has no poles for $|x|<1$.   Since $\rx(\omega;t)=F_{1}(x(\omega;t);t)$,
we find that $\rx(\omega;t)$ has no poles when $|x(\omega;t)|<1$.   With $\Lambda (\mathcal{O})=\{(x,y)\in \Etproj| |x|<1 \hbox{ or } |y|<1 \}$,  we deduce that a pole of $\rx(\omega;t)$  inside $\mathcal{O}$ satisfies $|y(\omega;t)|<1$.  We use 
 $\ry(\omega;t)=F_{2}(y(\omega;t);t)$,  and the fact that $F_{2}(y;t)$ has no poles for $|y|<1$
 to deduce that  $\ry(\omega;t)$ has no poles when $|y(\omega;t)|<1$. Therefore, the poles of 
 $\rx(\omega;t)$  inside $\mathcal{O}$ corresponds to the poles of $x(\omega;t)y(\omega;t)$ with $|y(\omega;t)|<1$.  The result follows. 
\end{proof}

With $\displaystyle \bigcup_{\ell\in \Z}\widetilde{\sigma}^{\ell}(\mathcal{O})=\C$ and $\widetilde{\sigma}^{-1}(\mathcal{O})\cap \mathcal{O}\neq \varnothing$, we then extend $\rx(\omega;t)$ and $\ry(\omega;t)$ as meromorphic functions on $\C$ where they satisfy the functional equations
\begin{align}
     \rx(\omega+\omega_{3}(t);t) & =\rx(\omega;t)+b_{x}(\omega;t) ,\label{eq:omega_3_per_rx} \\
     \rx(\omega+\omega_{1}(t);t)&=  \rx(\omega;t),\label{eq:omega_1_per_rx}    \\ 
     \ry(\omega+\omega_{3}(t);t)   &=\ry(\omega;t)+b_{y}(\omega;t), \label{eq:omega_3_per_ry} \\ 
     \ry(\omega+\omega_{1}(t);t)&=  \ry(\omega;t),\label{eq:omega_1_per_ry}
\end{align}
where $b_{x}(\omega;t) =y(-\omega;t)(x(\omega;t)-x(\omega+\omega_{3}(t);t))$ and $b_{y}(\omega;t)=x(\omega;t)(y(\omega;t)-y(-\omega;t))$.\smallskip

From  the functional equations~\eqref{eq:omega_1_per_rx} and~\eqref{eq:omega_1_per_ry}, the set of poles of $\omega\mapsto \rx(\omega;t)$ and $\omega\mapsto \ry(\omega;t)$ are $\omega_1(t)$ periodic.  With the  other functional equations 
and $\displaystyle \bigcup_{\ell\in \Z}\widetilde{\sigma}^{\ell}(\mathcal{O})=\C$, we may deduce the expression of a discrete set containing the poles of $\rx$ and $\ry$. 

\begin{lemma}\label{lem5}
Let $\mathcal{P}_{x}$ be the poles of $\rx$ in $\mathcal{O}$ and $\mathcal{P}_{b,x}$ be the poles of $b_{x}$ in $\C$. The set of poles of $\omega\mapsto \rx(\omega;t)$ is included in 
 $\left(\mathcal{P}_{x}+\omega_3 (t)\Z \right) \bigcup \left(\mathcal{P}_{b,x}+\omega_3 (t) \Z\right)$. A similar statement holds for $\ry (\omega;t)$.
\end{lemma}

\section{Preliminary results on differential algebraicity}\label{sec3}

In this section, we prove some results on differential algebraicity, and more specifically on $\partial_{t}$-algebraicity of the functions that appear in Section~\ref{sec2}.\smallskip

Let us begin by definitions. 
Let $f(x_1,\dots, x_n)$ be a  multivalued Puiseux series.  For $i=1,\dots,n$, we say that $f$ is $\partial_{x_{i}}$-algebraic if and only if it satisfies a nontrivial algebraic differential equation in the variable $x_i$,  with coefficients in $\Q$. We say that $f$ is differentially algebraic in all its variables (or differentially algebraic for short) if and only if for all $1\leq i\leq n$, $f$ is $\partial_{x_{i}}$-algebraic. \smallskip

The following remark, 
proved e.g.~in~\cite[Proposition 8, page 101]{Kol73},
will be used several times in the sequel.

\begin{remark}\label{rem2}
Let $f_1,\dots,f_n$ be differentially algebraic functions meromorphic on a common domain.
A function satisfies a nontrivial algebraic differential equation with coefficients in $\C(f_1,\dots,f_n)$ if and only if it satisfies a nontrivial algebraic differential equation with coefficients in $\Q$.
\end{remark}
The following lemma shows that the set of differentially algebraic functions is stable under many operations.

\begin{lemma}[Closure properties]\label{lem8}
The set of differentially algebraic functions meromorphic on a domain  is a field stable under   derivations.  If $f$ and $g$ are differentially algebraic and $f\circ g$ is well-defined then $f\circ g$ is  differentially algebraic as well. If $f$ is differentially algebraic  and admits an inverse $f^{-1}$, then $f^{-1}$ is also differentially algebraic.
\end{lemma}

\begin{proof}
See~\cite[Lemma 6.4]{DHRS} for the inverse property in the univariate case.  The proof extends straightforwardly to the multivariate case.  The rest of the statements follows from~\cite[Corollary 6.4 and Proposition 6.5]{bernardi2017counting}.
\end{proof}

In what follows, we might also consider functions of $t$ defined only on some intervals of~$(0,1)$.
Let ${\mathfrakD}$ be the field of  multivalued functions
that admit an expansion as convergent Puiseux series for all $t\in (0,1)$,  and that are differentially algebraic. In the sequel, when we will say that a function of $t$ defined (a priori) only of some intervals of $(0,1)$ is differentially algebraic, it will be implicit that it  may be continued as an element of ${\mathfrakD}$.\smallskip

The goal of the following results is to prove that various functions that appear in the uniformization of the elliptic curve  are $\partial_t$-algebraic.

\begin{lemma}[{\!\!\cite[Lemma~6.10]{bernardi2017counting}}]\label{lem6}
The functions $\omega_{1}(t),\omega_{2}(t),\omega_{3}(t)$ belong to~${\mathfrakD}$.\footnote{They even are solutions of linear differential equations.} Moreover, they  are analytic on a complex neighborhood of $(0,1)$.
\end{lemma}

\begin{proposition}\label{prop1}
Functions of ${\mathfrakD}(\wp(\omega;t),\partial_{\omega}\wp(\omega;t))$ are differentially algebraic in~$t$~and~$\omega$. 
\end{proposition}

\begin{proof}
Since the differentially algebraic functions form a field stable under the derivations (see Lemma~\ref{lem8}),  it suffices to show that $\wp (\omega;t)$ is differentially algebraic. 
It is well known that for $t\in (0,1)$ fixed, $\wp (\omega;t)$ is $\partial_{\omega}$-algebraic.  More precisely, it satisfies an equation of the form $(\partial_{\omega}\wp)^{2}=4\wp^{3}-g_{2}(t)\wp-g_{3}(t)$,
where $g_2(t),g_3(t)$  are the invariants of the elliptic curve.  Differentiating with respect of $\omega$ allows us to eliminate the invariants,  and obtain
$\partial_{\omega}^3 \wp (\omega;t)=12 \wp (\omega;t)\partial_{\omega} \wp (\omega;t)$; see~\cite[(18.6.5)]{abramowitz1988handbook}.
Hence $\wp (\omega;t)$ is $\partial_{\omega}$-algebraic.\par 
Let us prove the $\partial_t$-algebraicity.
In virtue of~\cite[Proposition 6.7]{bernardi2017counting}, $\wp$ satisfies a nontrivial $\partial_t$-algebraic equation with coefficients in 
$\C(\omega_1 (t),\omega_2 (t))$. By Lemma~\ref{lem6}, the periods $\omega_1 (t)$ and $\omega_2 (t)$ of $\wp$ are differentially algebraic, so in virtue of Remark~\ref{rem2}, $\wp (\omega;t)$ is $\partial_t$-algebraic.  
\end{proof}

\begin{remark}\label{rem3}
The same result holds 
with $\wp$ replaced by 
the Weierstrass function associated to the lattice $\omega_1(t)\Z+k\omega_2\Z$,
or the lattice $\omega_1(t)\Z+\omega_3 (t)\Z$.
\end{remark}

\begin{definition}[Principal part]
Let $ f(\omega;t)$ be a meromorphic function at $\omega=a(t)\in {\mathfrakD}$,  with Laurent series $f(\omega;t)=\sum_{\ell= \nu}^{\infty}a_{\ell}(t)(\omega-a(t))^{\ell}$.  The principal part of $f$ at~$\omega=a(t)$ is the sum $\sum_{\ell= \nu}^{-1}a_{\ell}(t)(\omega-a(t))^{\ell}$ with the convention that it is $0$ when $\nu \geq 0$.  The coefficients of this principal part are $a_{\nu}(t),\dots, a_{-1}(t)$.
\end{definition}
The following lemma will be used several times in the sequel. 

\begin{lemma}\label{lem7}
The following statements  hold:
\begin{itemize}
\item Let  $d(t)\in {\mathfrakD}$ be an arbitrary function. We have ${\wp(\omega;t)\in {\mathfrakD}((\omega+d(t)))}$;
\item $\wp(\omega;t)\in \omega^{-2}{\mathfrakD}[[\omega]]$;
\item The coefficients of the principal parts of $\omega\mapsto \wp (\omega;t)$ belong to~${\mathfrakD}$.
\end{itemize}
\end{lemma}

\begin{proof}
The last two assertions are straightforward consequences of the first one.  Let us prove the first point. 
The function $d(t)$ and the poles of $\omega\mapsto \wp(\omega;t)$   are analytic on a convenient domain. So either $-d(t)$ is a pole of $\omega\mapsto \wp(\omega;t)$ with constant order with respect to $t$, or the set of $t$ such that $-d(t)$ is a pole of $\omega\mapsto \wp(\omega;t)$ is discrete. It follows that the order of the pole of $\omega\mapsto \wp(\omega;t)$ at $-d(t)$ is constant except on a discrete set.
Since for $t$ fixed,  $\omega\mapsto \wp(\omega;t)$ has pole of order at most two,   we may write $\wp(\omega;t)=\sum_{\ell=k}^{\infty} c_{\ell}(t)(\omega+d(t))^{\ell}$. 

Note that the coefficients $c_{\ell}(t)$ may have a  pole when the order of the pole of $\omega\mapsto \wp(\omega;t)$ at $d(t)$ increases. 
In virtue of the field property of Lemma~\ref{lem8}, combined with  Proposition~\ref{prop1},  we find that $(\omega+d(t))^{-k}\wp(\omega;t)$  is differentially algebraic.   Note that
$c_{k}(t)$ is the value at $-d(t)$ of the $\partial_t$-algebraic function  $(\omega+d(t))^{-k}\wp(\omega;t)$.  By the field property of Lemma~\ref{lem8} and Proposition~\ref{prop1},  $(\omega+d(t))^{-k}\wp(\omega;t)$ is differentially algebraic in its two variables.   By the composition property of Lemma~\ref{lem8} it follows that $c_{k}(t)\in {\mathfrakD}$.
Let us fix $k \leq n$ and 	assume that  for $\ell=k,\dots ,n$, $c_{\ell}(t)\in {\mathfrakD}$. Let us show that $c_{n+1}(t)\in {\mathfrakD}$. This will prove the result by induction.
Let us define $h_{n}(\omega;t)=\wp(\omega,t)-\sum_{\ell=k}^{n}c_{\ell}(t)(\omega+d(t))^{\ell}$.
By Proposition~\ref{prop1},  the field property of Lemma~\ref{lem8}, and the induction hypothesis,  the function $t\mapsto h_{n}(\omega;t)$ is differentially algebraic in its two variables. Note that
$c_{n+1}(t)$ is the value at $-d(t)$ of $(\omega+d(t))^{-(n+1)}h_{n}(\omega;t)$.  By the composition property of Lemma~\ref{lem8} it follows that $c_{k}(t)\in {\mathfrakD}$.
\end{proof}

As a consequence of what precedes, we deduce:

\begin{corollary}\label{cor1}
The following holds.
\begin{itemize}
 \item The functions $x(\omega;t)$ and $y(\omega;t)$ are differentially algebraic in their two variables;
 \item For $d(t)\in {\mathfrakD}$, we have $x(\omega;t),y(\omega;t)\in {\mathfrakD}((\omega+d(t)))$;
 \item  The poles  and the coefficients of the principal parts of  $\omega\mapsto x(\omega;t)$ and $\omega\mapsto y(\omega;t)$ belong to~${\mathfrakD}$. 
\end{itemize}
\end{corollary}

\begin{proof}
We use the expressions of $x(\omega;t)$ and $y(\omega;t)$ given in Proposition~\ref{prop:uniformization}. The elements involved in the expression are meromorphic on some complex neighborhood of $(0,1)$ in the $t$-plane and are differentially algebraic by Proposition~\ref{prop1}. Since the differentially algebraic elements form a field, see Lemma~\ref{lem8}, the first point follows. Using Lemma~\ref{lem7}, we deduce that $x(\omega;t),y(\omega;t)\in {\mathfrakD}((\omega+d(t)))$ for all $d(t)\in {\mathfrakD}$.  Then the coefficients of the principal parts of  $\omega\mapsto x(\omega;t)$ and $\omega\mapsto y(\omega;t)$ belong to~${\mathfrakD}$. \par 
It remains to prove the differential algebraicity of the poles. 
Let $a(t)$ be a pole of  $\omega\mapsto x(\omega;t)$ or $\omega\mapsto y(\omega;t)$.  Then $a(t)$ is a  continuous function solution of $\wp(a(t);t)=b(t)$, where $b(t)$ is $\partial_t$-algebraic. 
\par 
Assume first that $\partial_\omega\wp(a(t);t)$ is identically zero or $a(t)$ is a pole of  $\wp(\omega;t)$.  By~\cite[page~270]{stein2010complex},  this corresponds to the case where 
$a(t)\in \omega_1(t)\frac{\Z}{2} +\omega_2(t)\frac{ \Z}{2} $.  By Lemma~\ref{lem6}, $a(t)$ is meromorphic on a complex neighborhood of $(0,1)$ and is $\partial_t$-algebraic. Then, it belongs to~${\mathfrakD}$. \par 
Assume now that $\partial_\omega\wp(a(t);t)$ is not identically zero and  $\wp(a(t);t)=b(t)$.
Then, by the implicit function theorem, $a(t)$ admits an expansion as a  meromorphic function on a complex neighborhood of any $t\in (0,1)$ with $\partial_\omega\wp(a(t);t)\neq 0$. On that domain $\wp$ is locally invertible and its inverse is differentially algebraic in its two variables by Lemma~\ref{lem8}. So we may write $\wp^{-1}(b(t);t)=a(t)$, where  $\wp^{-1}$ is the local inverse of $\wp$.
 With the composition and inverse properties of Lemma~\ref{lem8}, we deduce that $a(t)$ is $\partial_t$-algebraic.
Furthermore, by the implicit function theorem, it admits an expansion as a convergent series on a complex neighborhood of any $t\in (0,1)$.
 The  set of $t$ such that $\partial_\omega\wp(a(t);t)\neq 0$ being dense, we find that the differential equation holds everywhere. This concludes the proof.
\end{proof}
\pagebreak 

Recall, see Section~\ref{sec16}, that  $$b_{x}(\omega;t) =y(-\omega;t)(x(\omega;t)-x(\omega+\omega_{3}(t);t)) \hbox{  and }b_{y}(\omega;t)=x(\omega;t)(y(\omega;t)-y(-\omega;t)).$$ 
\begin{corollary}\label{cor2}
The following holds.
\begin{itemize}
 \item The functions $b_{x}(\omega;t)$ and $b_{y}(\omega;t)$ are differentially algebraic in their two variables;
 \item For $d(t)\in {\mathfrakD}$, we have $b_{x}(\omega;t),b_{y}(\omega;t)\in {\mathfrakD}((\omega+d(t)))$;
 \item  The poles  and the coefficients of the principal parts of  $\omega\mapsto b_{x}(\omega;t)$ and $\omega\mapsto b_{y}(\omega;t)$ belong to~${\mathfrakD}$. 
\end{itemize}
\end{corollary}

\begin{proof}
By Lemma~\ref{lem6},  $\omega_{3}(t)$ belongs to~${\mathfrakD}$.  This is now a straightforward application of Corollary~\ref{cor1}, combined with the field property of Lemma~\ref{lem8}. 
\end{proof}

Toward the proof of Theorem~\ref{thm1}, we are going to face to many situations where the series is known to be $\partial_{x}$-algebraic (or $\partial_{y}$-algebraic) for all fixed values $t$. More precisely the differential algebraicity of the series will be proved to be equivalent to the existence of functions that are for all $t$ fixed,  elliptic functions.  Unfortunately, few things are known about the $t$-dependence of the coefficients. 
The following result will be the main ingredient in the proof of Theorem~\ref{thm1} since it gives a framework where we can state that the elliptic functions are differentially algebraic in all their variables.  \par

\begin{theorem}\label{thm2}
Let $\omega \mapsto f(\omega;t)$ be a function such that:
\begin{itemize}
\item For all $t\in (0,1)$, $\omega\mapsto f(\omega;t)\in \C(\wp (\omega;t),\partial_{\omega}\wp (\omega;t))$. 
\item There are countably many elements of ${\mathfrakD}$, whose union forms the set of  poles  of $\omega\mapsto f(\omega;t)$. 
\item The coefficients of the principal parts  of $\omega\mapsto f(\omega;t)$ are in ${\mathfrakD}$. 
\item There exists $a(t)\in {\mathfrakD}$ such that  $f(a(t);t)\in {\mathfrakD}$. 
\end{itemize}
Then, $f(\omega;t)$ is differentially algebraic in its two variables.
\end{theorem}

\begin{remark}
At first sight, nothing is explicitly assumed on the $t$-dependence of $t \mapsto f(\omega;t)$. However, the assumptions on the poles, on the  principal parts,  and on the special value $f(a(t);t)$,  will imply that $t \mapsto f(\omega;t)$ is analytic on a convenient domain. 
\end{remark}
\begin{proof}
If $f$ is constant in the $\omega$ variable, then the result is clear. Assume that  $f(\omega;t)$ is not constant.  Let $a\in \C$.  By the field property in Lemma~\ref{lem8},   $f(\omega+a;t)$ satisfies the assumptions of Theorem~\ref{thm2}.  By the composition property,    $f(\omega;t)$ is differentially algebraic in its two variables if and only if  $f(\omega+a;t)$ is differentially algebraic in its two variables. Then, without loss of generality,  we may reduce to the case where for any pole $b(t)$ of $\omega \mapsto f(\omega;t)$,   $\partial_{\omega}\wp(b(t);t)$ is not identically zero.  We may also assume that $a(t)$ is not a pole of $\omega \mapsto \wp(\omega;t)$. 
\par 
Let us begin with the case where $\omega\mapsto f(\omega;t)$ is an even function.  As we can see in the proof of~\cite[Lemma 1.9]{stein2010complex}, we may write 
\begin{equation*}f(\omega;t)=c(t)\displaystyle \prod_{i=1}^{\kappa_{z}}f_{i}(\omega;t)\prod_{j=1}^{\kappa_{p}}g_{j}(\omega;t),
\end{equation*}\vspace{-5mm}

\noindent where  
\begin{itemize}
\item $c(t)$ is a function that does not depend upon $\omega$;
\item $f_{i}(\omega;t)$ are of the form $\wp(\omega;t)-\wp(a(t);t)$, where $a(t)$ are zeros of $\omega\mapsto  f(\omega;t)$;
\item $g_{j}(\omega;t)$ are of the form $(\wp(\omega;t)-\wp(b(t);t))^{-1}$, where $b(t)$ are poles  of $\omega\mapsto f(\omega;t)$.
\end{itemize}
\pagebreak

Then,  a partial fraction decomposition yields  a sum of the form 
\begin{equation}\label{eq6}
f(\omega;t)=\widetilde{c}(t)+ \sum_{i= 1}^{n_{\infty}}a_{i,\infty}(t)\wp (\omega;t)^{i}+ \sum_{j} \displaystyle \sum_{i= 1}^{n_j} \frac{a_{i,j}(t)}{(\wp (\omega;t)-\wp(b_j(t);t))^{i} }.
\end{equation}

By assumption,  the $b_j (t)$ are differentially algebraic.
Recall,  see Lemma~\ref{lem7},  that for all~$j$, we have ${\wp(\omega;t)\in {\mathfrakD}((\omega+b_j(t)))}$ (resp.~$\wp(\omega;t)\in \omega^{-2}{\mathfrakD}[[\omega]]$). Then, for every $i,j$, 
$$\frac{a_{i,j}(t)}{\Big(\wp (\omega;t)-\wp(b_j(t);t)\Big)^{i} }=\frac{a_{i,j}(t)}{\Big(\partial_{\omega}\wp(b_j(t);t)(\omega-b_j(t))\Big)^{i}}+O((\omega-b_j(t))^{-i+1}).$$

By the composition property of Lemma~\ref{lem8} and Proposition~\ref{prop1},  for all $k,\ell$ the function  $\left(\partial_{\omega}^{k}\wp(b_j(t);t)\right)^{\ell}$  is differentially algebraic.  Let us write the Taylor expansion of the function  \begin{equation*}
f(\omega;t)= \sum_{i=-n_j}^{\infty} \widetilde{a}_i (t)(\omega- b_j(t))^{i}.
\end{equation*}
Then,   for $i<0$, one has 
$$\widetilde{a}_i (t)=\frac{a_{i,j}(t)}{\partial_{\omega}\wp(b_j(t);t)^{i}}+f_{i,j},  \text{ where } f_{i,j}\in {\mathfrakD}( a_{i+1,j}(t),\dots ,a_{n_j ,j}(t)).$$   Since the coefficients of the principal part at $b_j(t)$ 
are differentially algebraic we have 
$\widetilde{a}_i (t)\in {\mathfrakD}$.
By  Lemma~\ref{lem8}, ${\mathfrakD}$ is a field,  and we find by a decreasing induction that for all $1\leq i \leq n_j$,   $a_{i ,j}(t)\in {\mathfrakD}$.     Similarly,  for all $i$,  we have 
$$  a_{i,\infty} (t)\wp (\omega;t)^{i}=\omega^{-2i} a_{i,\infty}(t)+O(\omega^{-2i+1}).$$  Then the coefficient of the term in  $\omega^{-2i}$ with $i>0$  is of the form $a_{i ,\infty}(t)+f_{i}$,  where $f_{i}\in {\mathfrakD}(a_{i+1,\infty }(t),\dots ,a_{n_\infty ,\infty}(t))$. 
  Since the coefficients of the principal part at $0$  are differentially algebraic,   we find $a_{i ,\infty}(t)+f_{i}\in {\mathfrakD}$.  By Lemma~\ref{lem8}, 
${\mathfrakD}$ is a field,  and  we find by a decreasing induction that for all $1\leq i \leq n_\infty$,   $a_{i ,\infty}(t)\in {\mathfrakD}$.  
Recall that by assumption, $f(a(t);t)$ is $\partial_t$-algebraic.  By Lemma~\ref{lem8} and Proposition~\ref{prop1},  we find \begin{equation*}\widetilde{d}(t):=\sum_{i= 1}^{n_{\infty}}a_{i,\infty}(t)\wp (a(t);t)^{i}+ \sum_{j} \displaystyle \sum_{i= 1}^{n_j} \frac{a_{i,j}(t)}{(\wp (a(t);t)-\wp(b_j(t);t))^{i} }\in {\mathfrakD}.
\end{equation*}
By the subtraction property of Lemma~\ref{lem8} we deduce that  $\widetilde{c}(t)=f(a(t);t)-\widetilde{d}(t)$ is $\partial_t$-algebraic.   In virtue  of Lemma~\ref{lem8} and Proposition~\ref{prop1},   every term in the right-hand side of~\eqref{eq6} is differentially algebraic.  With the field property of Lemma~\ref{lem8}, this concludes the proof in the  even case.  \smallskip

Assume that $\omega\mapsto f(\omega;t)$ is odd. 
The function $ \partial_{\omega}\wp(\omega;t)^{-1} f(\omega;t)$ is even, and  $\omega_1 (t)\Z+\omega_2 (t)\Z$, the poles of $\partial_{\omega}\wp(\omega;t)$, are $\partial_t$-algebraic; see Lemma~\ref{lem6}.  Then,  we may apply the even case to deduce that  
$ f(\omega;t)$ is of the form \begin{equation*}
\partial_{\omega}\wp(\omega;t)\widetilde{c}(t)+  \sum_{i= 1}^{n_{\infty}}a_{i,\infty}(t)\partial_{\omega}\wp(\omega;t) \wp (\omega;t)^{i}+ \sum_{j} \displaystyle \sum_{i= 1}^{n_j} \frac{a_{i,j}(t)\partial_{\omega}\wp(\omega;t)}{(\wp (\omega;t)-\wp(b_j(t);t))^{i} }.
\end{equation*}
\pagebreak 

Setting $a_{0,\infty}(t):=\widetilde{c}(t)$, we may rewrite the latter expression as
\begin{equation*}
\sum_{i= 0}^{n_{\infty}}a_{i,\infty}(t)\partial_{\omega}\wp(\omega;t) \wp (\omega;t)^{i}+ \sum_{j} \displaystyle \sum_{i= 1}^{n_j} \frac{a_{i,j}(t)\partial_{\omega}\wp(\omega;t)}{(\wp (\omega;t)-\wp(b_j(t);t))^{i} }.
\end{equation*}
\noindent By Proposition~\ref{prop1},  for all  $j$,  we have $\partial_{\omega}\wp(\omega;t)\in {\mathfrakD}((\omega-b_j(t)))$ (resp.~we have $\partial_{\omega}\wp(\omega;t)\in {\mathfrakD}((\omega))$).
The same reasoning as in the even case shows that  for all $1\leq i \leq n_j$,   the functions $a_{i ,j}(t)$  are differentially algebraic.  Similarly,   for all $0\leq i \leq n_\infty$, the functions $a_{i ,\infty}(t)$  are differentially algebraic.  By Proposition~\ref{prop1} and Lemma~\ref{lem8},  we find that 
 $f(\omega;t)$ is differentially algebraic.  This completes the proof in the odd case. \smallskip
 
Let us consider the general case. Note that by Proposition~\ref{prop1}, $\wp(\omega;t)-\wp(a(t);t)$ is differentially algebraic. So for all $n$, Lemma~\ref{lem8} ensures that $f(\omega;t)$ is differentially algebraic if and only if $(\wp(\omega;t)-\wp(a(t);t))^n  f(\omega;t)$ is differentially algebraic. So without loss of generality, we may reduce to the case where  $ f(\pm a(t);t)=0$.  We write 
$f(\omega;t)=f_{+}(\omega;t)+f_{-}(\omega;t)$, where 
\begin{align}
  f_{+}(\omega;t)&:=\frac{f(\omega;t)+f(-\omega;t)}{2},\\
  f_{-}(\omega;t)&:=\frac{f(\omega;t)-f(-\omega;t)}{2}.
\end{align}
  
The poles of $\omega \mapsto   f_{\pm}(\omega;t)$ are poles of $f$ or opposite of the latter. 
By Lemma~\ref{lem8}, they are $\partial_t$-algebraic  and  the coefficients of the principal parts are in ${\mathfrakD}$. Since  $f(\pm a(t);t)=0$ we find $f_{\pm}(a(t);t)=0$. In particular it is differentially algebraic. From the even and the odd cases, $f_{\pm}(\omega;t)$ are differentially algebraic in their two variables.
  Since the sum of two differentially algebraic functions is differentially algebraic, see Lemma~\ref{lem8},  we deduce that $f(\omega;t)=f_{+}(\omega;t)+f_{-}(\omega;t)$ is differentially algebraic.  
\end{proof}

\begin{remark}\label{rem4}\hspace{2em}
\begin{itemize}
\item As in Remark~\ref{rem3},  we may consider $\widetilde{\wp} (\omega;t)$,  the Weierstrass functions associated to  the lattice $\omega_1(t)\Z+\omega_3 (t)\Z$, or the lattice $\omega_1(t)\Z+k\omega_2 (t)\Z$, with $k\in \N^*$.  
Then, the proof of Theorem~\ref{thm2} can be straightforwardly adapted to this new lattice.  We then deduce the following. 
If  $\omega \mapsto f(\omega;t)$ is a function such that:
\begin{enumerate}
\item For all $t\in (0,1)$, $\omega\mapsto f(\omega;t)\in \C(\widetilde{\wp}(\omega;t),\partial_{\omega}\widetilde{\wp}(\omega;t))$. 
\item There are countably many elements of ${\mathfrakD}$, whose union forms the set of  poles  of  $\omega\mapsto f(\omega;t)$. 
\item The coefficients of the principal parts  of $\omega\mapsto f(\omega;t)$ are in ${\mathfrakD}$. 
\item There exists $a(t)\in {\mathfrakD}$ such that  $f(a(t);t)\in {\mathfrakD}$. 
\end{enumerate}
Then, $f(\omega;t)$ is differentially algebraic in its two variables.  
\item 
Let us now just assume that $\omega \mapsto f(\omega;t)$ satisfies the above  first three points and let $a(t)\in {\mathfrakD}$ that is not a pole.  Then,  $f(\omega;t)-f(a(t);t)$ satisfies the four points  and is therefore differentially algebraic.  By construction, the function $f(\omega;t)-f(a(t);t)$ has the same principal parts as $f(\omega;t)$.
\end{itemize}
\end{remark}

Although $\rx$ and $\ry$ are not elliptic functions, we will see in the next section that it is sufficient to control the behavior of their poles and coefficients in order to apply Theorem~\ref{thm2}.
\pagebreak 

\begin{lemma}\label{lem9}
The following holds:
\begin{enumerate}[label=\textnormal{(A{\arabic*})}]
\item The poles and 
coefficients of the principal parts of $\omega\mapsto \rx (\omega;t)$ belong to  ${\mathfrakD}$.\label{caseB} 
\item There exists $a(t)\in {\mathfrakD}$ such that  $ \rx (a(t);t)\in {\mathfrakD}$.\label{caseC} 
\end{enumerate}
Similar statements hold for $\ry$.
\end{lemma}

\begin{proof}
Let us prove the result for $\rx$, the reasoning for $\ry$ is similar. 
We refer to  Section~\ref{sec16} for the notations used in this proof.  

Recall that the series $Q(x,y;t)$ converges for $|x|,|y|,|t|<1$.  Let us consider $t$ in  $(0,1)$.  
Take $\omega\in\mathcal{O}$   (note that $\mathcal{O}$ depends continuously on $t$),
for each of the domains $|x(\omega;t)|<1$ and
$ |y(\omega;t)|<1$, 
one has the following equality of functions:
\begin{equation*}
F_{1}(x(\omega;t);t)=\rx (\omega ;t) \text{\qquad and \qquad} F_{2}(y(\omega;t);t)=\ry (\omega ;t),
\end{equation*}
with no poles on these domains.
Via the equality $0=\rx (\omega ;t)+\ry (\omega ;t)-K(0,0;t) Q(0,0;t)+x(\omega;t)y(\omega;t)$, and Lemma~\ref{lem4} on the inclusion of poles,  we deduce that the poles inside $\mathcal{O}$ of $\omega\mapsto \rx (\omega ;t)$ are the poles  inside $\mathcal{O}$ of 
$\omega\mapsto x(\omega;t)y(\omega;t)$ with $|y(\omega;t)|<1$.  What is more,  on that domain, $\omega\mapsto x(\omega;t)y(\omega;t)$  and   $\omega\mapsto \rx (\omega ;t)$  have the same principal parts.  
By Corollary~\ref{cor1},  the poles  of $\omega \mapsto \rx (\omega ;t)$ inside $\mathcal{O}$ are differentially algebraic.  Furthermore,  the corresponding principal parts have differentially algebraic coefficients.\smallskip 

Recall,  see~\eqref{eq:omega_3_per_rx},  that  $\rx(\omega+\omega_3(t);t) =\rx(\omega;t)+b_{x}(\omega;t)$.
By Corollary~\ref{cor2},  the poles and the coefficients of the principal parts  of $\omega\mapsto b_{x}(\omega;t)$ are differentially algebraic. 
By Lemma~\ref{lem6}, $\omega_3(t)$ is differentially algebraic.  Recall that 
  $\displaystyle \bigcup_{\ell\in \Z}\widetilde{\sigma}^{\ell}(\mathcal{O})=\C$.    
 With~\eqref{eq:omega_3_per_rx}   and what precedes,  we get assertion~\ref{caseB}.  \smallskip
 
  It remains to prove assertion~\ref{caseC}.  To lighten the notations we omit the dependence in $t$ in what follows.  Let us write $K(x,y;t)=\widetilde{B}_0 (y)+x\widetilde{B}_1 (y)+x^{2}\widetilde{B}_2 (y)$.
  Let~$\omega_0(t) \in \mathcal{O}$~such~that 
  \begin{equation*}
      y(\omega_0)=0 \text{\qquad and  \qquad} 
 x(\omega_0 )=\frac{-\widetilde{B}_1 (y(\omega_0))+\sqrt{\widetilde{B}_1 (y(\omega_0))^{2}-4\widetilde{B}_0 (y(\omega_0))\widetilde{B}_2 (y(\omega_0))}}{2\widetilde{B}_2 (y(\omega_0))}.
  \end{equation*}
The $y$-valuation of $\widetilde{B}_2 (y) $ being at most two, we consider the following subcases. 
\begin{itemize}
    \item If it is $0$ or $1$,  the valuation of the algebraic function $y\times \frac{-\widetilde{B}_1 (y) +\sqrt{\widetilde{B}_1 (y)^{2}-4\widetilde{B}_0 (y) \widetilde{B}_2 (y) } }{\widetilde{B}_2 (y)}$ is  nonnegative and then $\omega_0$ is not a pole of $x(\omega ;t)y(\omega ;t)$. \item If it is $2$,   then $4\widetilde{B}_0 (y) \widetilde{B}_2 (y)$
converges to $0$ when $y$ goes to $0$ and hence the same holds for 
$-\widetilde{B}_1 (y) +\sqrt{\widetilde{B}_1 (y)^{2}-4\widetilde{B}_0 (y) \widetilde{B}_2 (y) }$.  \end{itemize} 
We further find that $y\times \left(-\widetilde{B}_1 (y) +\sqrt{\widetilde{B}_1 (y)^{2}-4\widetilde{B}_0 (y) \widetilde{B}_2 (y) } )\right)\in O(y^2)$.  In that case,  we find that $\omega_0$ is not a pole of $x(\omega ;t)y(\omega ;t)$ either. 
 With $K(0,0;t) Q(0,0;t)=F_{2}(y(\omega_0;t);t)=\ry (\omega_0 ;t)$,  and $0=\rx (\omega ;t)+\ry (\omega ;t)-K(0,0;t) Q(0,0;t)+x(\omega;t)y(\omega;t)$, we then find  $0=\rx (\omega_0 ;t)+x(\omega_0 ;t)y(\omega_0;t)$.  It then suffices to show that  $x(\omega_0 ;t)y(\omega_0;t)$ is differentially algebraic.   With the expression of $y(\omega_0;t)$ in Proposition~\ref{prop:uniformization},  we find that $\omega_0$ is solution of an equation of the form $\wp(\omega_0;t)=b(t)$ with $b(t)\in  {\mathfrakD}$. With the same reasons as in the proof of Corollary~\ref{cor1},  we find that $\omega_0$ is differentially algebraic,  and $x(\omega_0 ;t),y(\omega_0;t)\in  {\mathfrakD}$. Then $ \rx (\omega_0 (t);t)\in {\mathfrakD}$. This concludes the proof. 
\end{proof}
\pagebreak

The following result relates the differential transcendence of $Q(x,y;t)$ and the differential transcendence of $\rx (\omega ;t)$ and $\ry(\omega ;t)$. 
\begin{proposition}\label{prop3}
The following statements are equivalent.
\begin{itemize}
\item The generating function $Q(x,y;t)$ is differentially algebraic in its three variables.
\item The series $F_{1}(x;t)$ and $F_{2}(y;t)$  are differentially algebraic in their two variables.
\item The meromorphic continuations $\rx (\omega ;t)$  and  $\ry (\omega ;t)$  are differentially algebraic in their two variables.
\end{itemize}
\end{proposition}
\begin{proof}
If $Q(x,y;t)$ is differentially algebraic then $Q(x,0;t)$ is differentially algebraic.  Since $K(x,0;t)$ is differentially algebraic, we use the ring property of Lemma~\ref{lem8} to deduce that 
$F_{1}(x;t)= K(x,0;t)Q(x,0;t) $ is differentially algebraic. (The reasoning is similar for the  differential algebraicity of   $F_{2}(y;t)$).
Conversely, if $F_{1}(x;t)$ and $F_{2}(y;t)$  are differentially algebraic, then, by evaluation, so is $Q(0,0;t)$.  As  the right-hand side of the expression in Lemma~\ref{lem1} is a sum and product of elements that are differentially algebraic,  
it is differentially algebraic 
(by the the field property in Lemma~\ref{lem8}). Therefore, $K(x,y;t)Q(x,y;t)$ is differentially algebraic.  Thus,
  $Q(x,y;t)$ is differentially algebraic.  So the first two points are equivalent.\smallskip

 Assume that the series $F_{1}(x;t)$ is  differentially algebraic in its two variables. Recall that $F_{1}(x(\omega;t);t)=\rx (\omega ;t)$  
 where $x(\omega;t)$ is differentially algebraic; see  Corollary~\ref{cor1}.
By composition of differentially algebraic functions,  see Lemma~\ref{lem8}, $\rx (\omega ;t)$ is differentially algebraic.
Conversely, on a domain where $x(\omega;t)$ is invertible, its inverse is also differentially algebraic; see Lemma~\ref{lem8}. We conclude similarly that if $\rx (\omega;t)$ is differentially algebraic then $F_{1}(x;t)$ is differentially algebraic. A similar reasoning holds for the $y$ variable and we find that $F_{2}(y;t)$  is differentially algebraic if and only if $\ry (\omega ;t)$  is differentially algebraic. This proves the equivalence between the last two points. 
\end{proof}

\section{Differential algebraicity of the generating function}\label{sec4}

The goal of this section is to prove Theorem~\ref{thm1} (the $\partial_x$, $\partial_y$, and $\partial_t$ differential algebraicity are equivalent). By Lemma~\ref{lem2}, the result holds for~all degenerate cases. By  Lemma~\ref{lem3} and Proposition~\ref{prop2}, it also holds when $\Etproj$ \textit{is not} an elliptic curve. So we now prove the case where $\Etproj$ is an elliptic curve. Let~$G$~be~the~group of the walk (see Definition~\ref{group}). Our proof handles separately the cases  $|G|<\infty$ and  $|G|=\infty$.

\subsection{Finite group case}

\begin{proposition}\label{prop4}
Let us consider a nondegenerate model of walks, assume that $\Etproj$ is an elliptic curve and $|G|<\infty$. Then, $Q(x,y;t)$ is $\partial_{x}$-algebraic, $\partial_{y}$-algebraic and $\partial_{t}$-algebraic.
\end{proposition}
\begin{proof}
By Proposition~\ref{prop3}, it suffices to show that $\rx (\omega;t)$  and $\ry (\omega;t)$ are differentially algebraic in their two variables.   Let us only consider $\rx (\omega;t)$, the proof for $\ry (\omega;t)$ is similar.
Recall that the $\omega_i(t)$ are continuous and that  $\omega_{3}(t)\in (0,\omega_{2}(t))$ (see Equation~\eqref{omega3}). Since  $|G|<\infty$ and $\widetilde{\sigma}(\omega)=\omega+\omega_3 (t)$, there exist $k,\ell\in \N^*$ with $\gcd (k,\ell)=1$ such that $\omega_3 (t)/\omega_{2}(t)=k/\ell$.  By~\eqref {eq:omega_3_per_rx},  we have  $ \rx(\omega+\omega_{3}(t);t)  =\rx(\omega;t)+b_{x}(\omega;t)$, where $b_{x}(\omega;t) =y(-\omega;t)(x(\omega;t)-x(\omega+\omega_{3}(t);t))$.  
 Let us recall some notations borrowed from the proof of~\cite[Theorem 4.1]{dreyfus2017differential}.
It is shown that we may write a decomposition of the form
\begin{equation}\label{eq4}
\rx (\omega;t)=\psi (\omega;t)+ \PHI (\omega;t)\phi (\omega;t). \end{equation}

More precisely, 
\begin{itemize}
\item $\PHI(\omega;t)= \displaystyle \sum_{i=0}^{\ell -1} b_{x}(\omega+ i\omega_3 (t);t)$;
\item $\phi (\omega;t)=\frac{\omega_{1}(t)}{2\mathbf{i}\pi}\zeta (\omega;t)-\frac{\omega}{\mathbf{i}\pi}\zeta (\omega_{1}(t)/2;t)$, where $\zeta$ is an opposite of the antiderivative of the Weierstrass function with periods $\omega_{1}(t)$ and $k\omega_{2}(t)$, that is  $$
\begin{array}{l}
\zeta (\omega;t)=\frac{1}{\omega}+\\
 \displaystyle\sum_{(\ell_{1},\ell_{2}) \in \Z^{2}\setminus \{(0,0)\}} \left(\frac{1}{\omega +\ell_{1}\omega_{1}(t)+\ell_{2}k\omega_{2}(t)} -\frac{1}{\ell_{1}\omega_{1}(t)+\ell_{2}k\omega_{2}(t)}+
\frac{\omega}{(\ell_{1}\omega_{1}(t)+\ell_{2}k\omega_{2}(t))^{2}}
\right);
\end{array}
$$
\item for all $t\in (0,1)$, the function  $\omega \mapsto \psi (\omega;t)$ is $(\omega_{1}(t),k\omega_{2}(t))$-periodic.
\end{itemize}
The idea is to prove successively  that $\PHI (\omega;t)$,  $\phi (\omega;t)$ and 
$ \psi (\omega;t)$ are differentially algebraic.  We will also  see them as functions of $\omega$ and 
study their poles and principal parts.
\begin{center}
Step 1: Study of  $\PHI (\omega;t)$. 
\end{center}

\begin{lemma}\label{lem11}
The following holds. 
\begin{itemize}
\item There are countably many elements of ${\mathfrakD}$,  whose union forms the set of  poles  of $\omega\mapsto \PHI (\omega;t)$. 
\item The coefficients of the principal parts  of $\omega\mapsto  \PHI (\omega;t)$ are in ${\mathfrakD}$. 
\item $\PHI$ is differentially algebraic in its two variables.  
\end{itemize}
\end{lemma}

\begin{proof}
Recall, see Lemma~\ref{lem6},  that $\omega_3(t)\in  {\mathfrakD}$. 
We may combine Corollary~\ref{cor2} 	and Lemma~\ref{lem8},  to deduce that  the poles and the coefficients of the principal parts of $\omega \mapsto \PHI(\omega;t)$ are $\partial_t$-algebraic.  Furthermore  by Lemma~\ref{lem8} and Proposition~\ref{prop1},  $\PHI$ is differentially algebraic in its two variables.
\end{proof}

\begin{center}
Step 2:  Study of $\phi (\omega;t)$. 
\end{center}
Before proving that $\phi (\omega;t)$ is differentially algebraic,  let us  study $\zeta (\omega;t)$.

\begin{lemma}\label{lem10}
The following holds. 
\begin{itemize}
\item There are countably many elements of ${\mathfrakD}$,  whose union forms the set of  poles  of $\omega\mapsto \zeta (\omega;t)$. 
\item The coefficients of the principal parts  of $\omega\mapsto \zeta (\omega;t)$ are in ${\mathfrakD}$. 
\item $\zeta$ is differentially algebraic in its two variables.  
\end{itemize}
\end{lemma}

\begin{proof}
In virtue of Lemma~\ref{lem6}, the periods $\omega_1 (t),\omega_2 (t)$ are differentially algebraic.  
Then, the poles and the coefficients of the principal parts of  $\omega \mapsto \zeta (\omega;t)$ belong to~${\mathfrakD}$. \smallskip
 
Let $\widetilde{\wp}$ be the Weierstrass function with periods $\omega_1 (t) ,k\omega_2 (t)$ and  let us write the classical differential equation
\begin{equation}\label{eq3}
\partial_{\omega}\widetilde{\wp} (\omega;t)^{2}=4\widetilde{\wp}(\omega;t)^{3}-\widetilde{g_{2}}(t)\widetilde{\wp}(\omega;t)-\widetilde{g_{3}}(t). 
\end{equation}
By Remark~\ref{rem3},  $\widetilde{\wp}(\omega;t)=-\partial_{\omega}\zeta (\omega;t)$ is differentially algebraic in its two variables. Then,  $\zeta (\omega;t)$ is $\partial_{\omega}$-algebraic too.  Let us prove the $\partial_t$-algebraicity. 
Let us differentiate~\eqref{eq3}  with respect to $\partial_{\omega}$ and simplify by $\partial_{\omega}\widetilde{\wp}(\omega;t)$, to find 
$$2\partial_{\omega}^{2}\widetilde{\wp} (\omega;t)=12\widetilde{\wp}(\omega;t)^{2}-\widetilde{g_{2}}(t). $$
By Lemma~\ref{lem8},  for all $k\geq 0$,  the derivatives $\partial_{\omega}^{k}\widetilde{\wp} (\omega;t)$ are $\partial_t$-algebraic.  Since the $\partial_t$-algebraic functions form a ring,  see Lemma~\ref{lem8}, we deduce that $\widetilde{g_2}(t)$ is $\partial_t$-algebraic.  Using again the ring property of  Lemma~\ref{lem8} in~\eqref{eq3}, we deduce that $\widetilde{g_3}(t)$ is $\partial_t$-algebraic too.  
We may see  the elliptic functions as functions of $\omega$ and $\widetilde{g_2},\widetilde{g_3}$.
By~\cite[(18.6.19)]{abramowitz1988handbook}, 
\begin{equation}\label{eq2}
(\widetilde{g_{2}}^{3}-27\widetilde{g_{3}}^{2})\partial_{\widetilde{g_3}}\widetilde{\wp}=( 3\widetilde{g_{2}} \zeta-\frac{9}{2} \widetilde{g_{3}}\omega)  \partial_{\omega}\widetilde{\wp}+6
\widetilde{g_{2}}\widetilde{\wp}^2
-9 \widetilde{g_{3}}\widetilde{\wp} -\widetilde{g_2}^2.
\end{equation}
 We have $\partial_{t}\widetilde{\wp}=\partial_t \widetilde{g_3} \partial_{\widetilde{g_3}}\widetilde{\wp}$.  If $\partial_t \widetilde{g_3}=0$ then $\widetilde{\wp}$ does not depend on $t$.  In particular,  its poles are independent of $t$,  which implies that the periods $\omega_1 (t)$ and $k\omega_2 (t)$ are independent of $t$.  Then,  $\zeta(\omega;t)$ is independent of $t$ and therefore $\partial_t$-algebraic.  
We similarly deal  with the case $\partial_t \widetilde{g_2}= 0$.
So let us consider the situation where both functions $\partial_t \widetilde{g_2}$ and $\partial_t \widetilde{g_3}$ are not identically zero.  
By the derivation property of Lemma~\ref{lem8},  we deduce that $ \partial_{t}\widetilde{\wp},\partial_t \widetilde{g_3} $ are $\partial_t$-algebraic.
Since the $\partial_t$-algebraic functions form a field, see  Lemma~\ref{lem8},  we then find that  $\partial_{\widetilde{g_3}}\widetilde{\wp}=\partial_{t}\widetilde{\wp} /\partial_t \widetilde{g_3} $ is $\partial_t$-algebraic.    Since $\partial_t \widetilde{g_2}\neq 0$,  we are in the situation where $\widetilde{g_2}$ is not identically zero.  
With~\eqref{eq2},  and the field property of Lemma~\ref{lem8},  we deduce that $\zeta(\omega;t) $ is $\partial_t$-algebraic.  This completes the proof of the lemma.
\end{proof}

\begin{lemma}\label{lem12}
The following holds. 
\begin{itemize}
\item There are countably many elements of ${\mathfrakD}$,  whose union forms the set of  poles  of $\omega\mapsto \phi (\omega;t)$. 
\item The coefficients of the principal parts  of $\omega\mapsto \phi (\omega;t)$ are in ${\mathfrakD}$. 
\item $\phi$ is differentially algebraic in its two variables.  
\end{itemize}
\end{lemma}

\begin{proof}
 In virtue of  Lemma~\ref{lem6}, the period $\omega_1 (t)$ is differentially algebraic.  
 By Lemma~\ref{lem8},  and Lemma~\ref{lem10},  we find that $\phi (\omega;t)$  is differentially algebraic in its two variables. Furthermore,   the poles and the coefficients of the  principal parts of $\omega \mapsto \phi (\omega;t)$ belong to~${\mathfrakD}$. \end{proof}
 
 \begin{center}
Step 3: Study of $\psi (\omega;t)$. 
\end{center}

Let us now study  $\psi (\omega;t)$.
By Lemma~\ref{lem9} there exists $a(t)\in {\mathfrakD}$ such that  $ \rx (a(t);t)\in {\mathfrakD}$.  Furthermore,   the  poles and the coefficients of the  principal parts of
 $\omega \mapsto \rx (\omega;t)$ are $\partial_t$-algebraic.
 \smallskip
 
With~\eqref{eq4}, Lemma~\ref{lem11} and Lemma~\ref{lem12},  we deduce that the poles of $\omega \mapsto \psi (\omega;t)$ are $\partial_t$-algebraic,  and  the coefficients of the principal parts are  $\partial_t$-algebraic.  
Recall that for all~$t$,  $\omega \mapsto \psi (\omega;t)$ is $(\omega_{1}(t),k\omega_{2}(t))$-periodic.
By Remark~\ref{rem4}, we may build $\omega \mapsto \psi_0 (\omega;t)$,  that is differentially algebraic and $(\omega_{1}(t),k\omega_{2}(t))$-periodic,  with same principal parts as $\omega \mapsto \psi (\omega;t)$.  We have 
\begin{equation}\label{eq10}
\rx (\omega;t)=\psi (\omega;t)-\psi_0 (\omega;t) + \PHI (\omega;t)\phi (\omega;t)+\psi_0 (\omega;t).
\end{equation}
Note that by construction  $\omega \mapsto \psi (\omega;t)-\psi_0 (\omega;t)$ has no poles.  Since $\omega \mapsto \rx (\omega;t)$  has no poles at $a(t)$,  we deduce with~\eqref{eq10}, that 
$\omega \mapsto\PHI (\omega;t)\phi (\omega;t)+\psi_0 (\omega;t)$ has no poles at $a(t)$.
Since $\PHI (\omega;t)\phi (\omega;t)+\psi_0 (\omega;t)$ is differentially algebraic (as the sum of two differentially algebraic functions,  see Lemma~\ref{lem8}),  with no poles at $a(t)$, we find that its evaluation at $a(t)$ is differentially algebraic.   Since $ \rx (a(t);t)\in {\mathfrakD}$ we use the ring property in  Lemma~\ref{lem8} to deduce that 
 $\psi (a(t);t)- \psi_0 (a(t);t)$ is differentially algebraic.  
 \pagebreak
 
 Then, 
 $\psi (\omega;t)- \psi_0 (\omega;t)$ satisfies the assumptions of Theorem~\ref{thm2} (with $\omega_2(t)$ replaced by $k\omega_2 (t)$) and we deduce that it is differentially algebraic by Remark~\ref{rem4}.  By Lemma~\ref{lem8},   and the differential algebraicity of  $\psi_0 (\omega;t)$, we deduce that  $\psi (\omega;t)$ is differentially algebraic. 
 
  \begin{center}
Step 4: Study of $\rx (\omega;t)$. 
\end{center}
 Let us now finish the proof of Proposition~\ref{prop4}. 
  Since 
 $\psi (\omega;t)$,   
 $\PHI (\omega;t)$, and $\phi (\omega;t)$ are  differentially algebraic in their two variables, we conclude that 
 $\rx(\omega;t)=\psi (\omega;t) + \PHI (\omega;t)\phi (\omega;t)$ 
 is differentially algebraic as the sum and product of differentially algebraic functions; see  Lemma~\ref{lem8}.
\end{proof}

\subsection{Infinite group case}
It now remains to handle the case where the group has  infinite  order.  So let us consider a nondegenerate model of walks and assume that $\Etproj$ is an elliptic curve and  $|G|=\infty$. The equivalence between the $\partial_{x}$-algebraicity and the $\partial_{y}$-algebraicity can be straightforwardly deduced in this weighted context from the proof of~\cite[Proposition~3.10]{DHRS}. Let us see that the $\partial_{t}$-algebraicity implies the $\partial_{x}$-algebraicity.
If $Q(x,y;t)$ is $\partial_t$-algebraic, then $Q(x,0;t)$ is $\partial_t$-algebraic.
By~\cite[Theorem 3.12]{dreyfus2019length}, if $Q(x,0;t)$ is $\partial_t$-algebraic, then it is $\partial_x$-algebraic. In virtue of Lemmas~\ref{lem1} and~\ref{lem8},  we find that if $Q(x,0;t)$ is $\partial_x$-algebraic, then $Q(x,y;t)$ is $\partial_x$-algebraic. 
So to prove Theorem~\ref{thm1}, it now suffices to show the following result.

\begin{theorem}
Let us consider a nondegenerate model of walks, assume that $\Etproj$ is an elliptic curve and  $|G|=\infty$.  If $Q(x,y;t)$ is $\partial_{x}$-algebraic, then it is $\partial_{t}$-algebraic.
\end{theorem}

\begin{proof} 
By  Proposition~\ref{prop3}, it suffices to show that $\rx(\omega;t)$ and  $\ry(\omega;t)$ are differentially algebraic.  Let us consider $\rx(\omega;t)$, the proof for $\ry(\omega;t)$  is similar. 
By Proposition~\ref{prop2}, for all $t\in (0,1)$ fixed, $\Etproj$ is an elliptic curve.  Let $G_t$ be the group $G$ specialized at~$t$. The order of the group $G_t$ may depend upon $t$.  However by~\cite[Proposition 2.6]{DHRS},  see also~\cite[Proposition~14]{kurkova2012functions},  which can be straightforwardly extended in the weighted framework,  the set of $t\in (0,1)$ such that $G_t$ has  infinite order is dense. 
By assumption, for such $t$ fixed, $x\mapsto F_{1}(x;t)$ is $\partial_{x}$-algebraic.
By~\cite[Theorem 3.8]{hardouin2020differentially}, for all such $t$ fixed there exists a $(\omega_1 (t),\omega_2 (t))$-periodic function $\widetilde{g}(\omega;t)$, such that 
\begin{equation}\label{eq5}
b_x(\omega;t)=\widetilde{g}(\omega+\omega_3 (t);t)-\widetilde{g}(\omega;t).
\end{equation}
 By~\cite[Proposition 3.9]{hardouin2020differentially},
 there exist $ g(x;t)\in \C(x,t)$ and  
$h(y;t)\in\C(y,t)$ such that $g(x(\omega;t);t)=\widetilde{g}(\omega;t)$ and for all $(x,y)\in \Etproj$, 
\begin{equation}\label{eq7}
xy=g(x;t)+h(y;t).
\end{equation}
Since $g(x(\omega;t);t)=\widetilde{g}(\omega;t)$, we use Corollary~\ref{cor1} to deduce that we may continue $\widetilde{g}(\omega;t)$ in the $t$ variable. 

\begin{center} Step 1: Study of  $\widetilde{g}(\omega;t)$.
\end{center}
\begin{lemma}
The following holds. 
\begin{itemize}
\item There are countably many elements of ${\mathfrakD}$,  whose union forms the set of  poles  of $\omega\mapsto \widetilde{g} (\omega;t)$. 
\item The coefficients of the principal parts  of $\omega\mapsto \widetilde{g} (\omega;t)$ are in ${\mathfrakD}$. 
\item $\widetilde{g}$ is differentially algebraic in its two variables.  
\end{itemize}
\end{lemma}

\begin{proof}
We claim that  the poles of $\omega\mapsto \widetilde{g}(\omega;t)$ are of the form $\omega_0(t)+\ell\omega_3 (t)$, where $\omega_0(t)$ is a pole of $\omega\mapsto b_{x}(\omega;t)$ and $\ell\in \Z$. To the contrary, assume that $a(t)$ is a pole that is not of this form. Then $a(t)-\omega_3 (t)$ is a pole of $\omega\mapsto \widetilde{g}(\omega+\omega_3 (t);t)$,  that is not a pole of $\omega\mapsto b_{x}(\omega;t)$.  With~\eqref{eq5}, we find that $a(t)-\omega_3 (t)$ is a pole of $\omega\mapsto \widetilde{g}(\omega;t)$.  We prove successively that for all $\ell \geq 0$, $a(t)-\ell \omega_3 (t)$ is a pole of $\omega\mapsto \widetilde{g}(\omega;t)$. Since $\widetilde{g}(\omega;t)$ is $(\omega_1 (t),\omega_2 (t))$-periodic, $a(t)-\omega_3 (t)\N+\omega_{1}(t)\Z+\omega_2 (t)\Z$ are poles of $\omega\mapsto \widetilde{g}(\omega;t)$.
Since $|G|=\infty$ and $\widetilde{\sigma}(\omega)=\omega+ \omega_3 (t)$, the sets $A_{\ell}:=\{a(t)-\ell\omega_3 (t)+\omega_{1}(t)\Z+\omega_2 (t)\Z\}$, with $\ell\in \N$, are disjoint.
Then, the set of poles of $\omega\mapsto \widetilde{g}(\omega;t)$  possesses an accumulation point which contradicts that the function is meromorphic. This proves the claim.\smallskip

By Corollary~\ref{cor2},  the poles of $\omega\mapsto b_{x}(\omega;t)$ are $\partial_{t}$-algebraic.  By Lemma~\ref{lem6},  $\omega_3 (t)$ is $\partial_t$-algebraic too.  
With the claim,   it follows that the poles of $\omega\mapsto \widetilde{g}(\omega;t)$ are $\partial_{t}$-algebraic.  By Corollary~\ref{cor1},  the coefficients of the principal parts of $\omega\mapsto x(\omega;t)$ are $\partial_t$-algebraic.  With $g(x(\omega;t);t)=\widetilde{g}(\omega;t)$,    and 
$g(x;t)\in \C(x,t)$,  we deduce that  the coefficients of the principal parts of $\omega\mapsto\widetilde{g}(\omega;t)$ are $\partial_t$-algebraic.  
Finally $\widetilde{g}(\omega;t)$ is differentially algebraic, as the composition of differentially algebraic functions; see Lemma~\ref{lem8}. 
\end{proof}

\begin{center}
Step 2:  Study of  $\widetilde{f}(\omega;t):=\rx(\omega;t)-\widetilde{g}(\omega;t)$. 
\end{center}
By~\eqref{eq:omega_3_per_rx} and~\eqref{eq5}, we find 
\begin{align*}
\widetilde{f}(\omega+\omega_3(t);t)&=\rx(\omega+\omega_3(t);t)-\widetilde{g}(\omega+\omega_3(t);t)\\
&=\rx(\omega;t)+ b_{x}(\omega;t)-(\widetilde{g}(\omega;t)+b_{x}(\omega;t))=\widetilde{f}(\omega;t).
\end{align*}
Then, 
 $\widetilde{f}(\omega;t)$ is $\omega_3(t)$-periodic. Recall that $\widetilde{g}(\omega;t)$ is $\omega_1(t)$-periodic.
By~\eqref{eq:omega_1_per_rx},  $\rx(\omega;t)$ is also $\omega_1(t)$-periodic.  Therefore, $\omega\mapsto \widetilde{f}(\omega;t)$ is elliptic with periods $(\omega_1(t),\omega_3 (t))$.
Recall that  the poles and the coefficients of the principal parts of $\omega \mapsto  \widetilde{g}(\omega;t)$  are $\partial_t$-algebraic.  By Lemma~\ref{lem9}, the same holds for $\rx(\omega;t)$ and  there exists $a(t)\in  {\mathfrakD}$ such that $\rx(a(t);t)$ is differentially algebraic.  By Remark~\ref{rem4}, we may build $\omega \mapsto\widetilde{f}_0 (\omega;t)$,  that is differentially algebraic,   $(\omega_{1}(t),\omega_{3}(t))$-periodic,  and with same principal parts as $\omega \mapsto \widetilde{f} (\omega;t)$.  Let us write 
$$\widetilde{f}(\omega;t)- \widetilde{f}_0 (\omega;t) =\rx(\omega;t)-\widetilde{g}(\omega;t)- \widetilde{f}_0 (\omega;t).$$
The function $-\widetilde{g}(\omega;t)- \widetilde{f}_0 (\omega;t)$ is differentially algebraic, as the sum of two differentially algebraic functions; see Lemma~\ref{lem8}.  Since $\widetilde{f}(\omega;t)- \widetilde{f}_0 (\omega;t) $ has no poles and $a(t)$ is not a pole of 
 $\omega \mapsto \rx(\omega;t)$, we deduce that $a(t)$ is not a pole of  $\omega \mapsto -\widetilde{g}(\omega;t)- \widetilde{f}_0 (\omega;t)$. 
Therefore, its evaluation at $a(t)$ is still differentially algebraic.   
Then,   the same holds for $\omega \mapsto \widetilde{f}(\omega;t)- \widetilde{f}_0 (\omega;t)$, 
which satisfies the assumptions of Theorem~\ref{thm2} (with $\omega_2(t)$ replaced with $\omega_3 (t)$) and we deduce that it is differentially algebraic by Remark~\ref{rem4}.
Hence,  $\rx(\omega;t)=\left(\widetilde{f}(\omega;t)- \widetilde{f}_0 (\omega;t)\right)+\widetilde{g}(\omega;t)+ \widetilde{f}_0 (\omega;t)$ is differentially algebraic  as the sum of differentially algebraic functions,  see  Lemma~\ref{lem8}.   This concludes the proof.
\end{proof}

\subsection*{Acknowledgments}

This project has received funding from the ANR de rerum natura ANR-19-CE40-0018.
\pagebreak

\bibliographystyle{SLC}
\bibliography{Dreyfus}
\end{document}